\documentclass{article}

\usepackage[utf8]{inputenc}
\usepackage[english]{babel}

\usepackage{authblk}
\usepackage{graphicx}
\usepackage{amsmath,amssymb,amsfonts}
\usepackage{amsthm}
\usepackage{mathrsfs}
\usepackage[title]{appendix}
\usepackage{xcolor}
\usepackage{textcomp}
\usepackage{booktabs}
\usepackage{algorithm}
\usepackage{hyperref}
\usepackage{listings}
\usepackage{orcidlink}
\usepackage{eurosym}
\usepackage{url}
\usepackage{siunitx}

\graphicspath{
    {figures/}
}
\DeclareSIUnit{\sieuro}{\text{\euro}}

\title{Energy-Aware Aggregation of Input Data for the Optimisation of Heat Supply of Municipal Districts}

\author[1]{
    Patrik Schönfeldt~\orcidlink{0000-0002-4311-2753}
}

\author[1]{
    Elif Turhan~\orcidlink{0000-0001-9903-7628}
}

\affil[1]{German Aerospace Center (DLR), Institute of Networked Energy Systems, Oldenburg, Germany}

\begin{document}

\maketitle

\section*{Abstract}
In the context of municipal heat planning, it is imperative to consider the numerous buildings, numbering in the hundreds or thousands, that are involved.
This poses particular challenges for model-based energy system optimization, as the number of variables increases with the number of buildings under consideration.
In the worst case, the computational complexity of the models experiences an exponential increase with the number of variables.
Furthermore, within the context of heat transition, it is often necessary to map extended periods of time (i.e., the service life of systems) with high resolution (particularly in the case of load peaks that occur at the onset of the day).
In response to these challenges, the aggregation of input data is a common practice.
In general, building blocks or other geographical and urban formations, such as neighbourhoods, are combined.
This article explores the potential of incorporating energy performance indicators into the grouping of buildings.
The case study utilizes authentic data from the Neu-Schwachhausen district,
grouped based on geographical location, building geometry, and energy performance indicators.
The selection of energy indicators includes the annual heat consumption as well as the potential for solar energy generation. To this end, a methodology is hereby presented that considers not only the anticipated annual energy quantity, but also its progression over time.
We present a full workflow from geodata to a set of techno-socio-economically Pareto-optimal
heat supply options.
Our findings suggest that it is beneficial to find a balance between geographical position
and energy properties when grouping buildings for the use in energy system models.

\section{Introduction}

In Germany, the heating sector accounts for more than 84 \% of the final energy consumption of the residential sector and contributes about a third of the country’s GHG emissions~\cite{iea2020energy}.
To guide the necessary process of decarbonisation,
the heat planning act (``Wärmeplanungsgesetz'') requires every municipality to submit a heat plan~\cite{arnold2024warmewende}.
This type of heat planning needs to consider hundreds to thousands of buildings.
In particular for model based optimisation, this imposes a challenges,
as worst case the computational time can increase exponentially with the number of variables~\cite{Goldfarb1994}.
Even worse, in the context of heat planning, both long (component lifetimes) and short (peak loads, fluctuations of renewable generation) time scales are important.
To face these challenges, aggregation of data is required.
On the geographical scale, often building blocks,
neighbourhoods, or other groups are introduced by city planers~\cite{AMARAL2018406, BLANCO2024105075}.
Secondly, there are building typologies focusing solely on energy~\cite{GANGOLELLS20163895,LOGA20164}.
These typically focus on the building demand,
independent from potential characteristics of the building stock in a specific area (with exception for localized classification, e.g.~\cite{DALLO2012211}) or generation potential.
When it comes to generation potential,
some works investigate the potential solar gains depending on building types in a neighbourhood~\cite{10.1108/SASBE-08-2017-0035,su14159023}.
This approach, however, is more linked to a green field approach.

Several years ago, data acquisition could be considered the main issue
to elaborate heat transition plans for municipalities.
At that time, ``common data-driven models apply statistical and artificial intelligence (AI) simulation mainly based on machine learning (ML) methods to model energy use''~\cite{ABBASABADI2019106270}.
With the aforementioned law in action,
municipalities are now mandated to collect data for the creation of the heat plans.
Additionally, solar maps~\cite{KANTERS20141597} nowadays are widely implemented. 
Thus, the focus shifts to availability of planning capacities and on data processing,
facilitating a more data-centred approach.
From a general perspective,
the idea described in this contribution can be described as follows:
\begin{enumerate}\itemsep=0.1em
    \item compress the input data (cf. Sec.~\ref{sec: clustering}),
    \item process the compressed data (cf. Sec.~\ref{sec: system optimisation}),
    \item decompress result data to have result for original dataset (cf. Sec.~\ref{sec: data analysis}).
\end{enumerate}
If compression and processing the data are commutable,
the final result will be the same as if the result of the original data was compressed.
Although we are dealing with linear energy system models,
even a linear compression method will only fulfil this requirement
in edge cases where the energy system model does not reach any constraint.

In the following, we walk through the method using the area of Bremen
Neu-Schwachhausen and the corresponding energy geodata as an example.
Where meaningful to evaluate the method,
we applied simplifying assumptions.

\section{Characterisation of input data}

\subsection{Heat demand}
\begin{figure}[ht]
    \centering
    \includegraphics[width=0.7\linewidth]{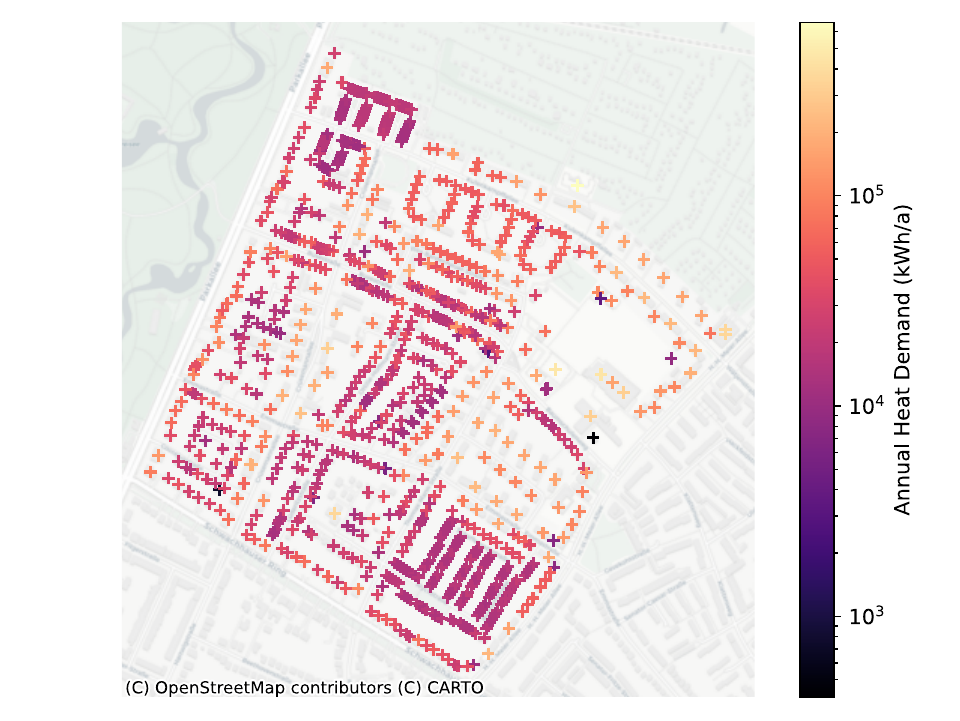}
    \caption{
        Annual heat demand of all 859 buildings in the area.
    }
    \label{fig: heat demand map}
\end{figure}

\begin{figure}[ht]
    \centering
    \includegraphics[width=0.48\linewidth]{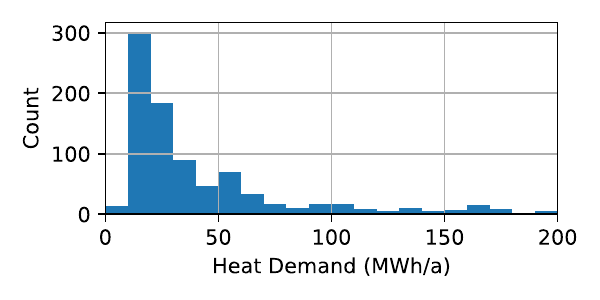}\hfill
    \includegraphics[width=0.48\linewidth]{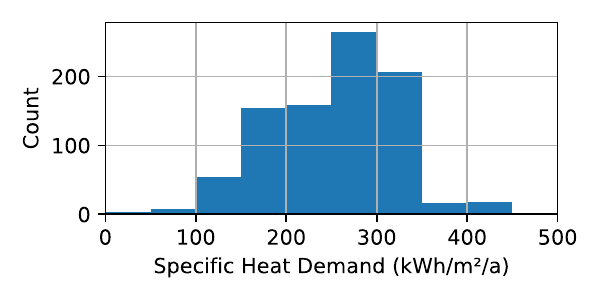}
    \caption{
        Histograms of the annual heat demand of the buildings
        in absolute numbers per building (left) and
        relative to the building footprint (right).
    }
    \label{fig: heat demand histograms}
\end{figure}

A map showing the annual heat demand of the individual buildings
is displayed in Fig.~\ref{fig: heat demand map}.
By the bare eye, it is visible that there are darker
(less heat demand) and lighter (higher heat demand)
areas.
This mostly reflects areas with bigger and smaller
building sizes and not the achieved building standard.
The influence of the building footprint is also dominating
the distribution of the annual heat demand
as is displayed in Fig.~\ref{fig: heat demand histograms}.
The curve for the absolute heat demand peaks before
\SI{20}{MWh/a}, then falls quickly with a long tail
starting at \SI{70}{MWh/a}, which continues to show non-zero
values even beyond the scale.
The complementary numbers relative to the building footprint can be seen in
the right panel of the same figure.
This distribution is very compact, outliers are rare,
the middle four of the the equally sized bins are all populated
by between 150 and 250 buildings.

\subsection{Solar potential}

When it comes to local energy generation,
solar potential marks a key parameter.
It can be calculated from roof parameters,
that are available, e.g. from LOD2 open data \cite{LOD2_Germany2024}.
For this study, we use data from the solar map of Bremen,
that also includes shadowing by trees.
As buildings in urban areas are aligned with streets,
which follow an almost orthogonal pattern in the region
we investigate,
the roof azimuth groups in four different orientations.
\begin{figure}[ht]
    \centering
   \includegraphics[width=0.48\linewidth]{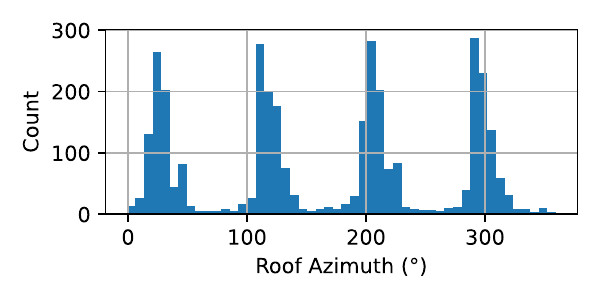}
    \hfill
   \includegraphics[width=0.48\linewidth]{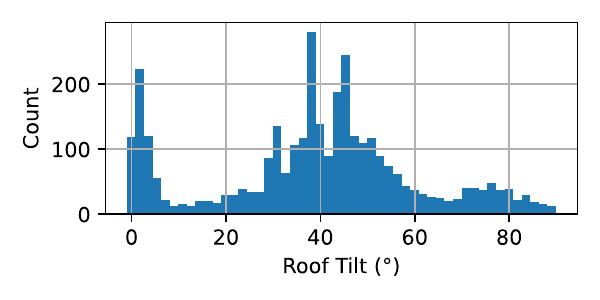}
    \caption{
        Histograms of the roof azimuth (left)
        and roof tilt (right).
    }
    \label{fig: roof histograms}
\end{figure}
As displayed in Fig.~\ref{fig: roof histograms},
these corresponding groups are rather distinct.
However, it has to be kept in mind that later not the individual
roofs but the buildings they cover are important.
The roof tilt follows a not so clear distribution,
an interpretation might need to include knowledge about
the psychology of architects.

However, what we are actually interested in is not the shape of the building
but its energy-related properties.
Thus, we apply a method to our knowledge first described in~\cite{dlr212738}.
The idea is to approximate the generation potential curve of every roof
by a linear combination of the generation potential curves of \(N\) `standard roofs'
\begin{equation}
    P_{\alpha, \theta}(t) \approx \sum_{n=0}^{N-1} w_{n} \times P_n (t).
\end{equation}
This way, there is a distinct quantity (the weight \(w_n\) of the `standard roof')
that clustering can use.
When compared to using the annual solar potential instead,
this new method preserves information about the temporal structure
on a finer temporal resolution.
Later, the factors \(w_n\) can be used with the solar potential time series \(P_n\)
to reproduce the solar potential time series for every building.
It should be mentioned that it is still an ad-hoc model,
that has been shown plausible but is missing a proper validation
according to scientific standards.

\begin{figure}[bth]
    \centering
    \includegraphics[width=0.7\linewidth]{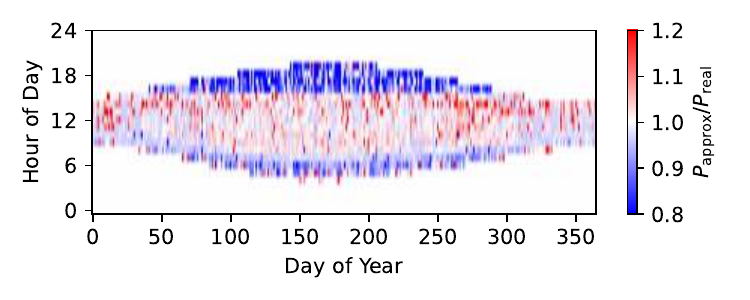}
    \caption{
        Fraction of approximated over real solar power potential values over the year.
    }
    \label{fig: pv waterfall}
\end{figure}

While it is possible to find the coefficients \(w_n\)
e.g. using a least square fit,
\cite{dlr212738} suggests to just solve a set of linear equations.
One argument to do so is the reduction of computational load,
which can be important especially for large datasets.
When picking a discrete number of times of the day
\(\SI{0}{h} \le t_\mathrm{d} < \SI{24}{h}\) and the same number of
representative roofs (defined by a combination of \(\alpha\) and \(\theta\)),
the set of equations to solve in order to obtain the \(w_n\) is defined by
\begin{subequations}
\begin{equation}
    \langle P_{\alpha, \theta}\rangle(t_\mathrm{d}) = \sum_{n=0}^{N-1} w_{n} \times \langle P_n \rangle (t_\mathrm{d})
    \quad \forall\; t_\mathrm{d}, P_{\alpha, \theta}, P_n
\end{equation}
where
\begin{equation}
    \langle P_x\rangle (t_\mathrm{d}) = \frac{1}{365} \times \sum_{d=0}^{364} P_x (d \times \SI{24}{h} + t_\mathrm{d}),
\end{equation}
with \(P_x \in \{P_{\alpha, \theta}, P_n\}\)
is the power at that time of the day averaged over the whole year.
\label{eq: representative roof weights}
\end{subequations}
(Note that Eq.~\ref{eq: representative roof weights} assumes consistent time throughout the year,
i.e. daylight saving time has to be considered manually.)
For the present study, the standard roofs are chosen to have a common tilt
of \(\theta = \SI{30}{\degree}\) and azimuth angles of \(\alpha \in (\SI{90}{\degree}, \SI{180}{\degree}, \SI{270}{\degree})\).
Further we choose \(t_\mathrm{d} \in (\SI{9}{h}, \SI{12}{h}, \SI{15}{h})\).
Results produced using this method are displayed in Fig.~\ref{fig: pv waterfall}.

For the present example, the approximation works quite well in the range defined
by the chosen \(t_\mathrm{d}\).
Here, the deviation is typically in the single per-cent range.
However, there is a systematic underestimation of the power
in early mornings and late evenings, particularly in summer.
As absolute numbers approach low values at these times,
we assume an acceptable impact of this bias.

\begin{figure}[ht]
    \centering
   \includegraphics[width=0.48\linewidth]{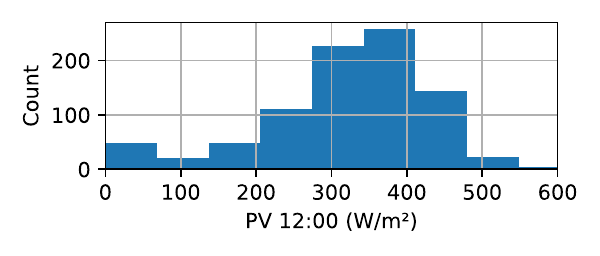}
    \caption{
        Histogram of the solar energy generation potential
        (relative to building footprint, not to roof area).
    }
    \label{fig: pv histogram}
\end{figure}

The distribution of a corresponding, derived quantity (PV potential at noon)
for all buildings in the area can be seen in Fig.~\ref{fig: pv histogram}.
This quantity, again, is more compact compared to the original quantities.
In our dataset, we also find four buildings with extraordinary high solar
potential relative to the building footprint
(\(P_\mathrm{peak} > \SI{600}{W/m^2}\)).
These values are outliers but no obvious error,
as the roof area can exceed the building footprint because of the roof tilt.
As the method we apply later is rather robust against these types of data errors,
we can just accept them.
For this particular area, we also observe that there is a strong linear
correlation between the solar potential at different times of the day.
As this might be due to the high amount of parallel buildings in the dataset
(cf. upper histogram in Fig.~\ref{fig: roof histograms}),
we do not generalise this finding and use
\(\langle P_{\alpha, \theta}\rangle(t_\mathrm{d})\) for all of the chosen
 \(t_\mathrm{d}\) to describe the solar potential.

\subsection{Cost of heat network connection}

The cost of heat network connection is per building,
but as we work with values relative to the building footprint,
we need to treat the costs for heat network connections the same way.
We estimate the costs of the heat network connection by first computing the
peak heat demand according to a BDEW (`Bundesverband der Energie- und Wasserwirtschaft', German for `National Association of Energy and Water Industries') demand profile assuming a single family
house~\cite{demandlib_v022}
\begin{equation}
    \hat{P}_\mathrm{th} [\si{kW}] \approx 0.228 \times E [\si{MWh/a}].
\end{equation}
This number is then used for a cost estimation of a heat network connection according to the
KWW Technikkatalog~\cite{kww_technikkatalog} as
\begin{equation}
        C_\mathrm{HNC}(\hat{P}_\mathrm{th}) =  \SI{13.4111}{\sieuro/kW} \times \hat{P}_\mathrm{th} + \SI{13976}{\sieuro}.
\end{equation}

\begin{figure}[ht]
    \centering
   \includegraphics[width=0.48\linewidth]{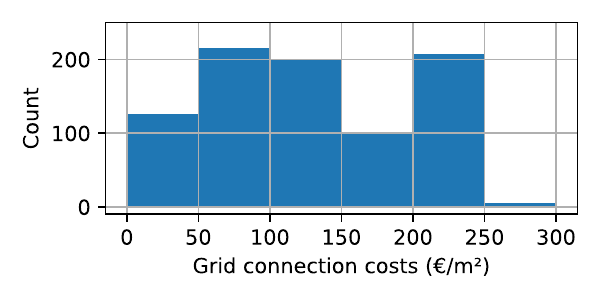}
    \caption{
        Histogram of the grid connection costs
        (relative to building footprint).
    }
    \label{fig: heat grid histogram}
\end{figure}

The distribution is displayed in Fig.~\ref{fig: heat grid histogram}.
There are four outliers, with values above \SI{300}{\sieuro/m^2}
up to close to \SI{600}{\sieuro/m^2}.
That this number equals the number of outliers with very high solar potential,
is pure coincidence: These high-cost buildings have a solar potential
of \SI{0}{kW/m^2}.

\section{Vector quantisation}
\label{sec: clustering}

We search for buildings that are similar,
i.e. that can meaningfully rely on the same heating solutions.
The similarity includes two aspects: location and energy.
Because clusters exist in \(n\)-dimensional data space,
they often overlap in their projection to a two-dimensional plane.
When focusing on decentralised heat supply options as in~\cite{dlr212738},
this is no issue.
However, it would be a bad assumption that supply using a new heating grid is equally suitable for every single building in a geographically scattered cluster.
In fact, we identify different requirements for the clustering algorithm for energy and spatial indicators:
Energy-wise, clusters should be as homogenous as possible to be able to propose the same building energy system for all of them.
Location-wise, the extent of the clusters might not be as crucial as long as the area is connected, i.e. suitable for a common centralised solution.
For all the methods, the data has been normalised to the range [0:1] using the \texttt{MinMaxScaler} of~\cite{scikit-learn}.

\subsection{Data-based building categories}

As a first step, we cluster based on energy values.
The goal is to automatically derive building categories
tailored to both the data and the optimisation process.
Because we want every building within one cluster to be represented
by one representative building, clusters should be compact
in the sense that the standard deviation within the cluster
should be low for every value.
This property is served by K-Means clustering.
As a side-note, clustering including integer dimensions is also possible when
K-Medoids~\cite{k_medoids_python} is used instead of K-Means.
This might be required, e.g. if the model included a discrete number
of inhabitants as a variable.

\begin{figure}[hbt]
    \centering
    \includegraphics[width=0.5\linewidth]{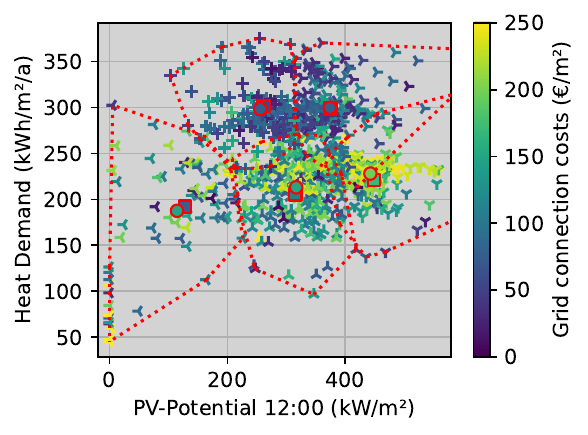}
    \caption{
        Scatter-plot of all buildings with respective footprint-specific
        energetic quantities.
        Five clusters are identified using different markers and grouped by lines.
        Possible cluster representatives are drawn as
        disks (average values per building, original K-Means centres)
        and as squares (building values weights by building footprint).
    }
    \label{fig: energy clusters}
\end{figure}

For the current case study, however, the considered values are positive real numbers.
We have specific heat demand in \si{kWh/m^2/a},
solar potential in the morning (9 am), at noon (12 pm), and in the afternoon (3pm),
all in \si{W/m^2},
and costs for a heat grid connection in \si{\sieuro/m^2},
where the \si{m^2} reflects the building footprint for all five quantities.

The result of an example clustering (five clusters)
using all of the aforementioned energetic properties 
is displayed in Fig.~\ref{fig: energy clusters}.
Note that the clusters overlap only in the 2D projection.
There are two high-demand clusters,
that both also feature relatively low grid connection costs,
one with medium and one with high solar potential.
These are complemented by three low-demand clusters
with higher grid connection costs, one of those also features
low solar potential.
That the cluster centre of the values relative to the building footprint
is very close to the building footprint-weighted average
corresponds to the fact that for this particular dataset
\begin{equation}
    \sum_{b\in\mathfrak{B}} A_b \times \sum_{b\in\mathfrak{B}} \frac{x_b}{A_b}
    \approx 
    \mathrm{N}_\mathrm{b} \times \sum_{b\in\mathfrak{B}} x_b,
\end{equation}
where \(\mathfrak{B}\) is the set of building indexes,
\(x_b\) is any quantity of building \(b\) (e.g. its annual heat demand), and
\(A_b\) is the footprint of that building.
This means that in this case the method using representatives expressed in specific quantities does not bias the summed values of the clusters.

\subsection{Geographical grouping}
\label{sec: geographical grouping}

In this second step,
buildings are grouped considering (also) their geographical locations.
These groups will be used to decide on heat-network connections.
As connections are not realised in a star topology but using line strips, 
compactness is not as important as for the building categories,
which leaves a rich set of options for potential algorithms,
including DBSCAN~\cite{10.5555/3001460.3001507},
OPTICS~\cite{10.1145/304182.304187}, and
HDBSCAN \cite{10.1007/978-3-642-37456-2_14}.
While we do not do this here,
note that these also allow for using a custom definition for the distance between two points.

While certain configurations -- particularly
doing the same thing for all buildings --
are possible independent from the grouping,
other options depend on the clustering method.
We did a scan over different numbers for both,
representative buildings and geographical groups.
This scan ranged from 5 to 30 for both,
representative buildings and geographical groups.
Results are summarized in Table~\ref{tab: 2-step clustering scan}.
A visualisation and deeper discussion of the results can be found in Appendix~\ref{sec: two-step scan}.
We find that there is a trade-off between the total number of variables needed for the model and the geographical cohesion
of the geographical groups.
We measure the latter by calculating the average length to connect buildings of only one geographical group to a heat network.
In the following,
they described further interwoven with the descriptions of the methods.

\begin{table}
    \centering
    \begin{tabular}{c|cc|cc}
        Method & \multicolumn{2}{c}{Variables} & \multicolumn{2}{c}{Network Length (\si{m})}\\
            &min    &max    &min    &max\\
        \hline
        K-Means         & 25    &   378 & 30.8& 31.9\\
        K-Means Energy  & 19    &	187 & 41.6&	64.5\\
        K-Prototypes    & 23    &	197 & 40.0&	88.1\\
        K-Prototypes after HDBSCAN
                        & 20    &	193 & 39.4&	80.9\\
        K-Modes         & 8     &	67  & 50.2&	119.3\\
    \end{tabular}
    \caption{
        Extreme values in the 2-step clustering scan
    }
    \label{tab: 2-step clustering scan}
\end{table}

\paragraph{K-Means} Clustering using geographic location. This leads to geographically dense, non-overlapping groups.
While the single-group line length is almost constant using this method,
increasing the number of representative buildings or the number of geographical groups leads to an almost proportional increase of the number of variables.

\paragraph{K-Means Energy} K-Means clustering using geographic location and energy data. These groups are still geographically dense but allow for some overlap to reflect similarities of buildings.
In our scan, this approach about halves the number of variables
at the cost that the average line length now increases with the number
of representative buildings.

\paragraph{K-Modes}
K-Modes~\cite{huang1998extensions},
generalises K-Means for categorical data.
In a first sub-step,
connected areas within buildings of the same building category are identified using HDBSCAN.
In a second sub-step,
the labels from the HDBSCAN runs and the building categories
(from step 1) are used to find the final groups.
As this method focuses on finding groups containing as few
representative buildings as possible,
it clearly reduces the number of variables in the model.
This comes at the cost of increased line length if only one
geographical group is connected to using heat network.

\paragraph{K-Prototypes} This algorithm combines K-Modes
with K-Means to handle datasets with mixed scalar and categorical data~\cite{huang1998extensions,devos2015}.
This way, building categories and geographical locations can be considered.
In the scan, an increase in the number of geographical clusters
has a bigger impact on the number of variables compared to
an increase of the number of representative buildings.
Line length also increases for both but not as much as with
K-Modes.

\paragraph{K-Prototypes after HDBSCAN} A combination of geographic location,
building categories, and labels of HDBSCAN is used.
Compared to pure {K-Prototypes},
this method comes with less increased intra-cluster
line length when both energy and final clusters are high.
On the other hand, the number of variables increases slightly
more with the number of representative buildings
(starting from a lower minimum).

\begin{figure}[thb]
    \centering
    \includegraphics[width=0.49\linewidth]{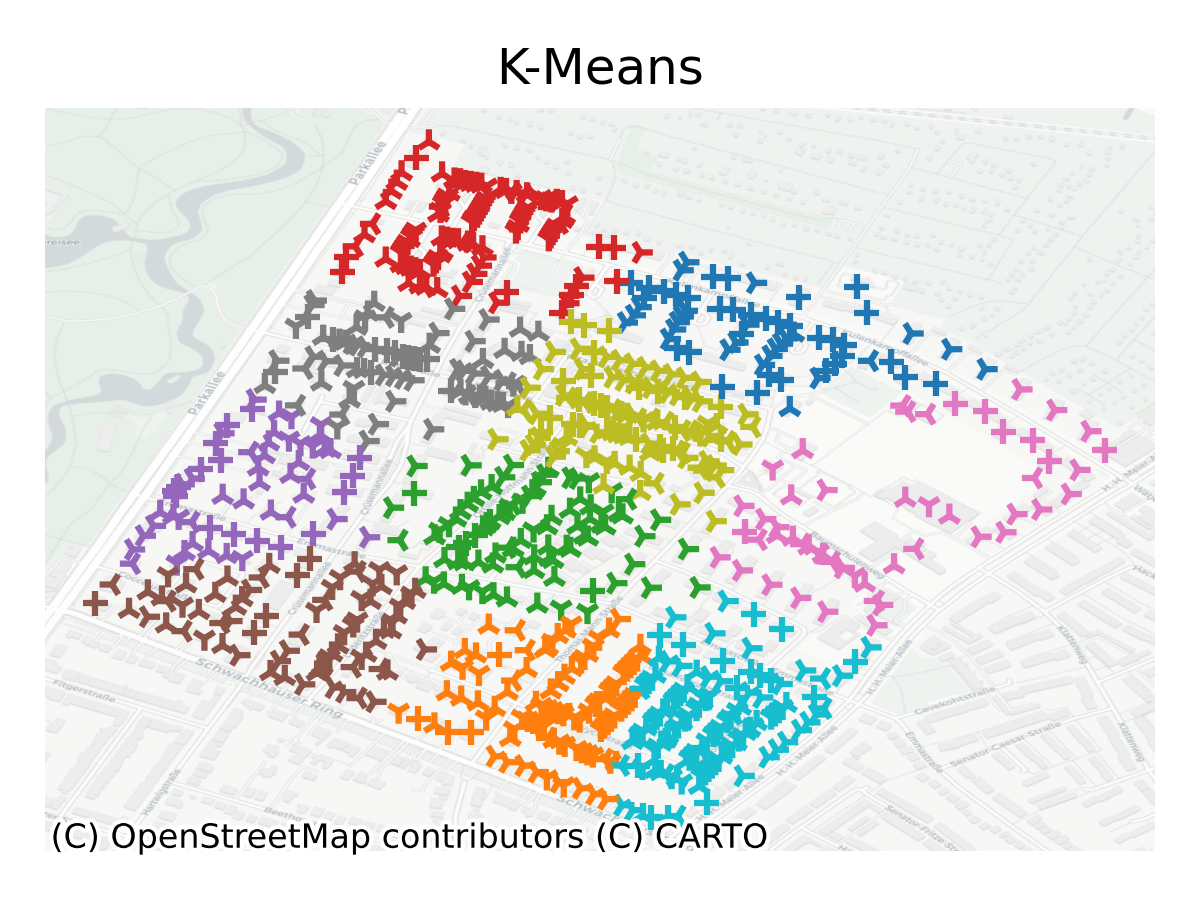}
    \includegraphics[width=0.49\linewidth]{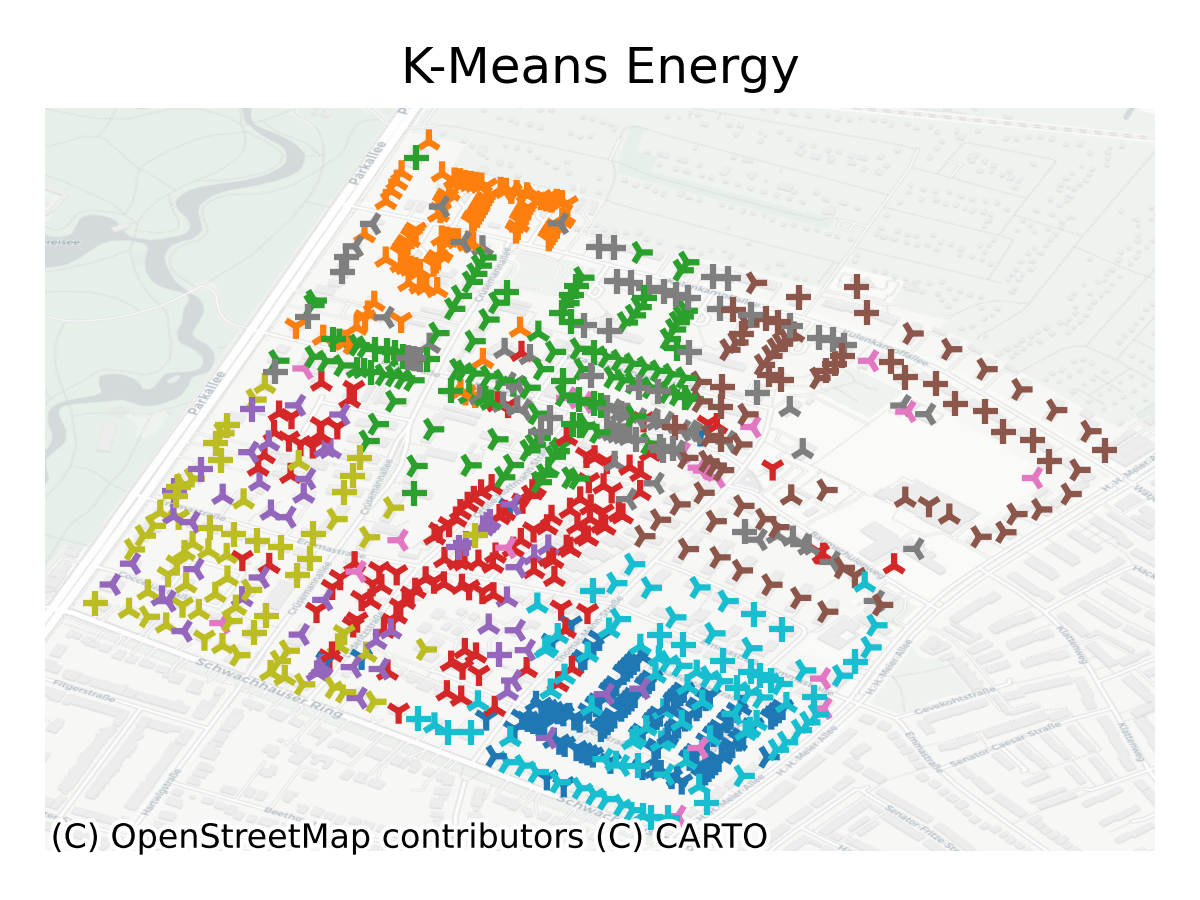}\\
    \includegraphics[width=0.49\linewidth]{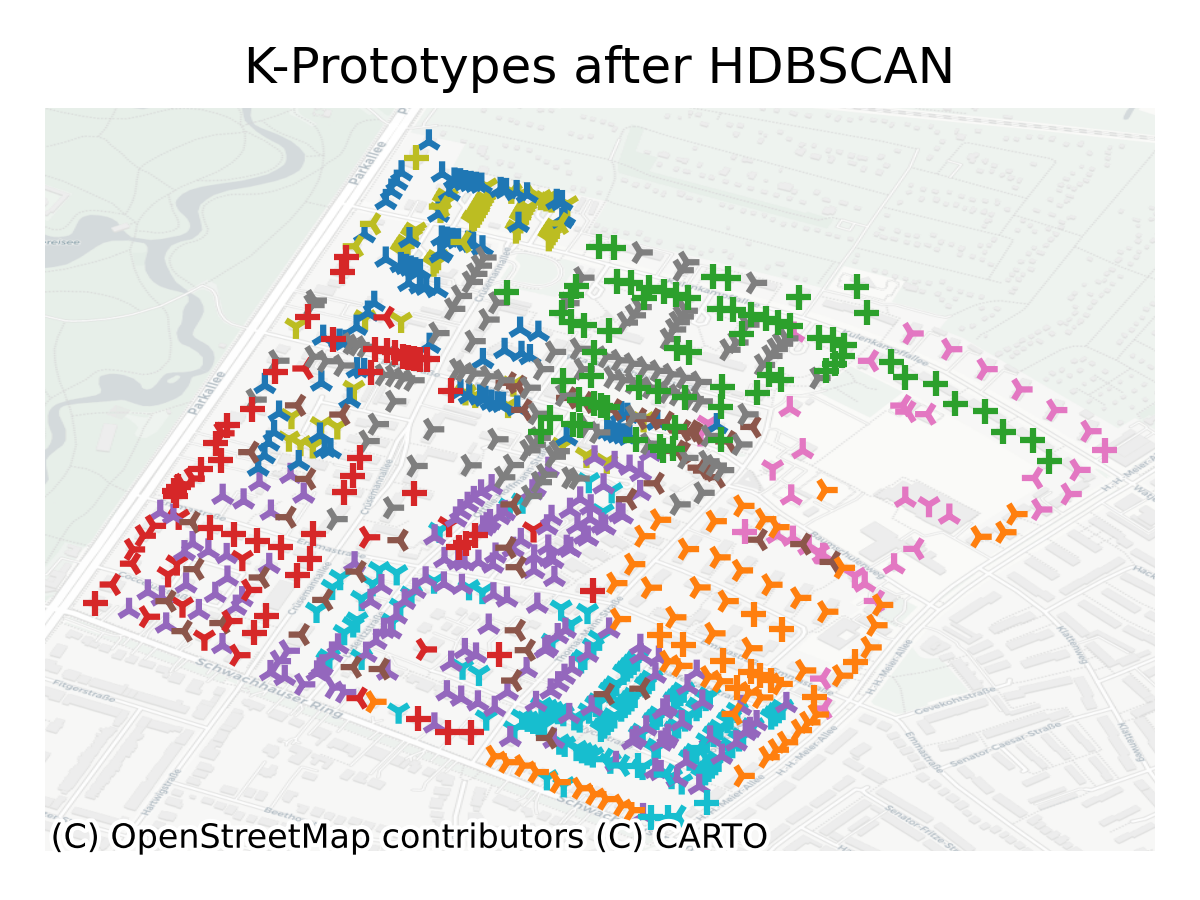}
    \includegraphics[width=0.49\linewidth]{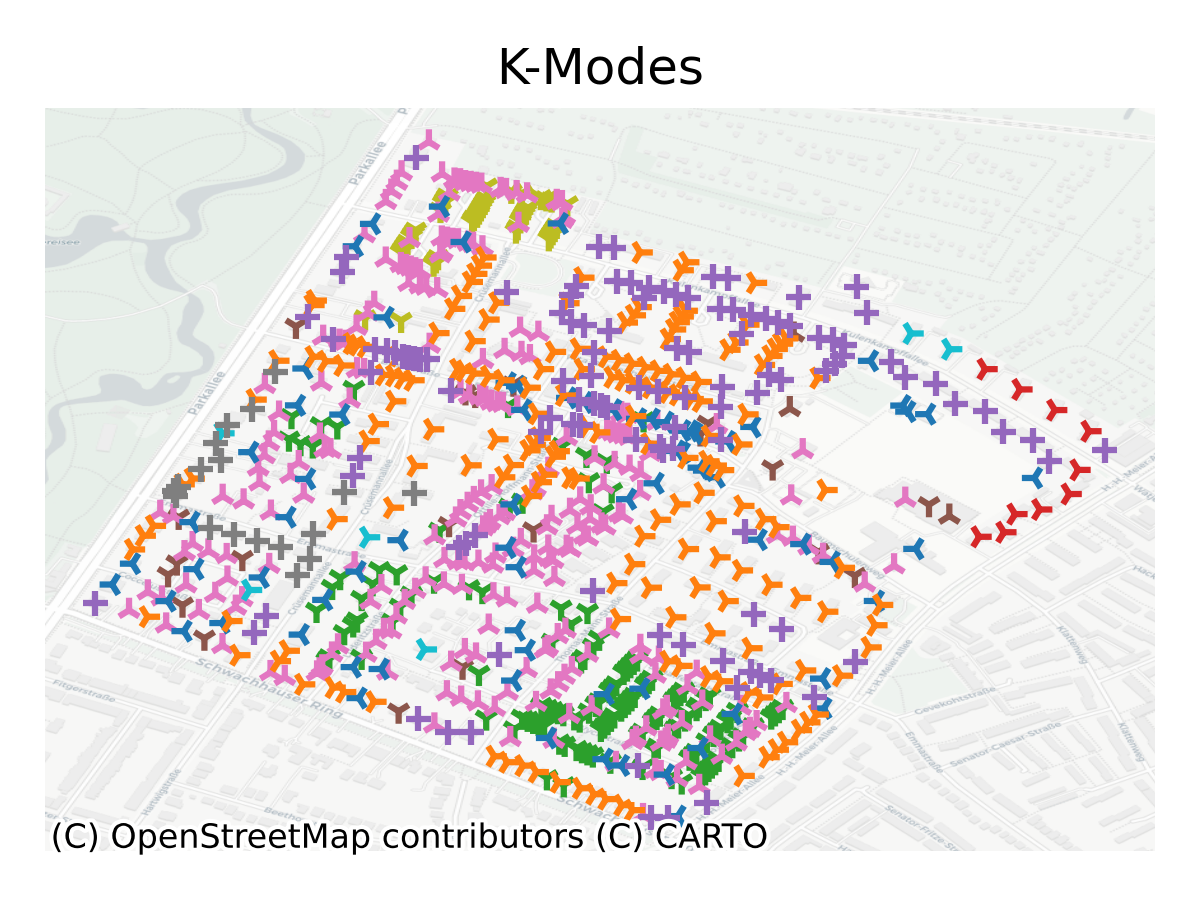}
    \caption{
        Assignment of buildings to representative buildings (icons)
        and geographical groups (colours) for the case of five
        representative buildings and ten geographical groups.
    }
    \label{fig: map clusters}
\end{figure}

As they provide very similar performance,
in the following, we drop the application of K-Prototypes
without the extra HDBSCAN.
An example assignment of the presented two-step methods
using five representative buildings and ten geographical groups
is shown in Fig.~\ref{fig: map clusters}.
It is visible that the methods have different geographical cohesion.
Still, the assignment is also clearly considering the location,
as seen in a number of very localised groups
(light green, red, grey) even in the K-Modes example.

\section{Optimisation model}
\label{sec: system optimisation}

The optimisation is an iterative process performed using a heuristic search.
For every archetype building,
we investigate the two options of a centralised (i.e. heat network) or a decentralised heat supply.
For both of them, the size of the PV system is optimised,
for the latter also the choice between installing a new air source heat pump
and keeping an (assumed to be existing) gas boiler.
If better data about currently existing heating technologies was at hand,
it should be considered in the clustering and would translate to price differences
for different technologies.

Note that we decided on parameters in a way that the influence of
the clustering method can be studied.
Thus, technologies should have clear advantages and disadvantages.
For example, keeping a gas boiler results in operation costs but no investment costs.
Secondly, there is a decision between the
centralised an the decentralised option of a building per geographical group.
This is where the number of combinations becomes relevant.
For the present study with five representative buildings
in ten geographical groups,
the number of combinations of representative buildings and geographical groups,
which is the number of decision variables in the optimisation,
as well as the shortest possible heat grid when selecting only one of these combinations
is shown in Tab.~\ref{tab: 10-5 variables}.

\begin{table}
    \centering
    \begin{tabular}{l c r}
        Method & Variables & Shortest grid\\
        \hline
        K-Means & 50 & \si{275}{m}\\
        K-Means Energy & 31 & \si{26}{m}\\
        K-Prototypes after HDBSCAN & 29 & \si{83}{m}\\
        K-Modes & 12 & \si{1754}{m}\\
    \end{tabular}
    \caption{
        Number of variables,
        as well as shortest possible heat grid
        (one group, one representative building).
    }
    \label{tab: 10-5 variables}
\end{table}

For the supply of the heat network,
we assume a central heat pump that has access to waste heat (COP of 8),
while decentralised heat pumps use the ambient air
(COP ranging between 2.4 and 6,
depending on the season and time of the day).
After the decision is made, which buildings are connected,
we use a minimum spanning tree~\cite{spanning_tree}
to find the length of the heat network.
To estimate the power required for groups of buildings using
a standard load profile, we apply an empirical formula for the diversity factor of \(n\) buildings~\cite{Winter2001-Gleichzeitigkeit},
\begin{equation}
    q_\mathrm{diversity}(n) = \min\left(
        1, a + \frac{b}{1 + (n/c)^d}
    \right),
\end{equation}
with
\(a = 0.450\),
\(b = 0.551\),
\(c = 53.8\), and
\(d = 1.76\) (rounding to three significant digits by us).
The formula originally targeted \( 1 < n \le 200\),
but as it approaches almost constant values afterwards,
we also use it for higher numbers of buildings.
In particular, by dividing the peak load of the standard load profile in the limit \(q_\mathrm{diversity}(n \rightarrow \infty) \approx 0.47\), we find the peak load for one building.
Note that the formula was derived using data from only two heat networks,
so the method has to be understood as a rough approximation.
The presented aggregation method works independently
from the chosen energy system model for the representative buildings.
Thus, we apply a rather simple model which is presented in
Appendix~\ref{sec: energy system model}.

\begin{subequations}
Using the number of the occurrences of the representative buildings
in the geographical groups \(N_{b,g}\), the global indicators
can be computed
for total energy costs in 2025
\begin{equation}
    C_{\mathrm{E},2025}
    = C_\mathrm{el,HN} + \sum_{b\in\mathfrak{B}}\sum_{g \in \mathfrak{G}} N_{b,g} C_{b,g}
    \label{eq: indicator energy costs 2025}
\end{equation}
with the set of group indexes \(\mathfrak{G}\),
the electricity costs for the heat network \(C_\mathrm{el,HN}\)
and the energy costs per building \(C_{b,g}\) (see Eq.~\ref{eq: building enery costs}),
and for annualised investment costs
\begin{equation}
    C_{\mathrm{invest}} = C_{\mathrm{HN}} + \sum_{b \in \mathfrak{B}} \sum_{g \in \mathfrak{G}} \sum_{i \in \mathfrak{I}} N_{b,g} C_{i,b,g},
    \label{eq: indicator investment costs}
\end{equation}
with the set of technologies to be invested in de-centrally \(\mathfrak{I}\),
their annualised costs \(C_{b,g,i}\).
At last, with \(P_{\mathrm{el}}(t)\) the electrical and \(P_{\mathrm{gas}}\)
the gas consumption time series if the complete area,
the green house gas emissions are calculated as
\begin{equation}
    m_{\mathrm{GHG}} = \int \frac{f_\mathrm{el}}{4} \max\left(
            P_{\mathrm{el}}(t), \SI{0}{kW}
        \right)
        + f_\mathrm{gas} P_{\mathrm{gas}}(t)\;\mathrm{d}t,
    \label{eq: indicator emissions}
\end{equation}
\end{subequations}
where \(f_\mathrm{el} = \SI{260}{g/kWh}\)
and \(f_\mathrm{gas} = \SI{240}{g/kWh}\)
are the specific emission factors for 2025 and the factor \(1/4\)
reflects the assumption that the average emissions will be
\(1/4\) of that value over the lifetime.

After all indicators are evaluated for the system,
the heuristic tries new values for the variables.
Details about this method can be found in~\cite{SCHMELING20221223}.
(As we decided to not include energy storage in the present study,
the operation is inflexible and not optimised here.)
For the present study,
we allow 100\;000 generations in a population of 16
using the Non-dominated Sorting Genetic Algorithm II~\cite{996017}
as available in PyGMO~\cite{Biscani2020}.
Leaving more time for the algorithm will increase
the likelihood of sampling more of the Pareto-plane.
However, we use the quality of the results as an
indicator to access the sampling methods.

The result of this method is a set of (nearly)
Pareto-optimal configurations for the energy-system.
It should be noted that the sampling of the solution space is random.
This implies that the density of points does not contain any information.
In particular,
a) the non-existence of a configuration in the result data does
not mean that it is not Pareto-optimal and
b) points at the edge of the sampled Pareto-hyperplane might not be
the most extreme possible options.
We will see a concrete example in the next section.

\section{Data analysis}
\label{sec: data analysis}

To get an overview over the multitude of Pareto-optimal solutions,
it makes sense to also work with aggregated values
instead of only looking at every single option for every single representative building in each of the geographical groups.
After the general overview,
we will look into example configurations
resulting from the optimisation.
These configurations can be seen as the decompressed results,
as they view at every individual building again.

\begin{figure}[thb]
    \centering
    \includegraphics[width=0.49\linewidth]{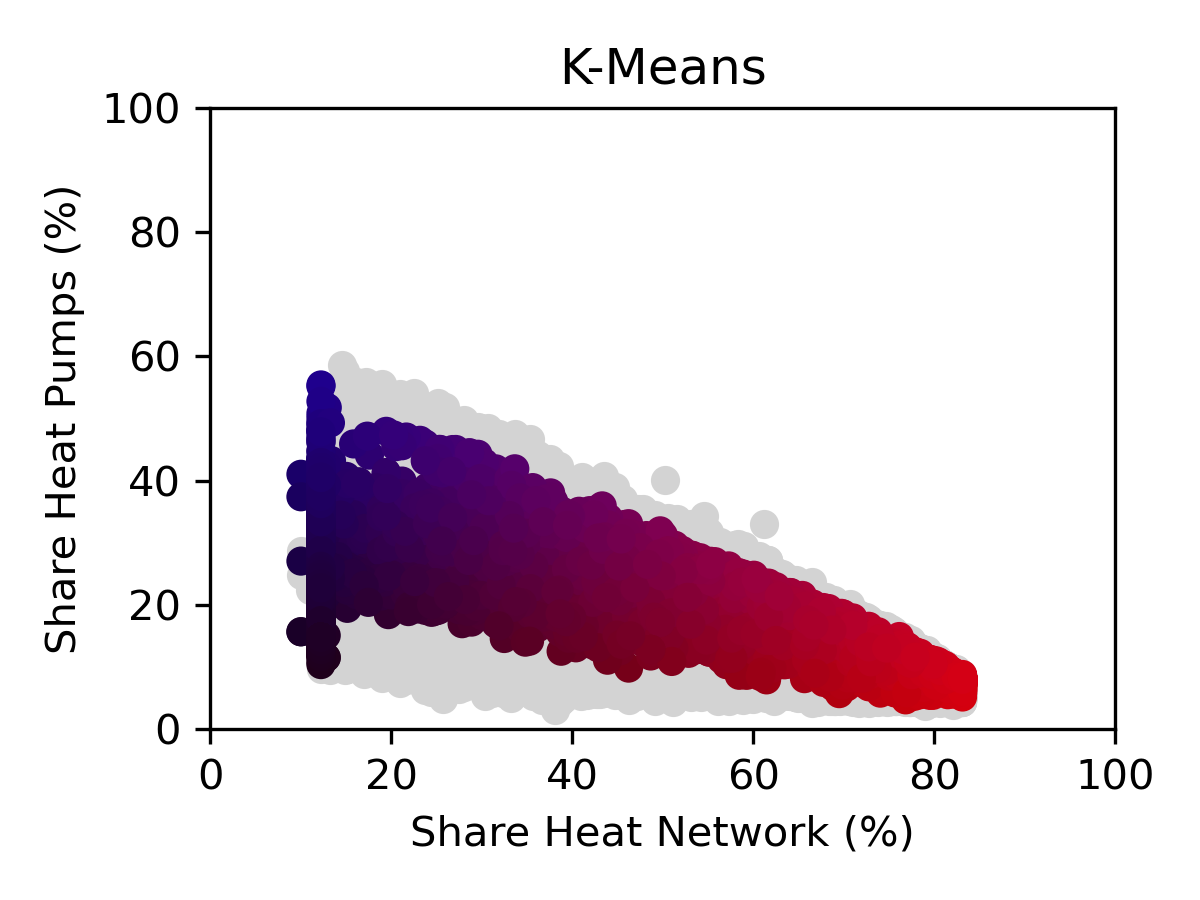}
    \includegraphics[width=0.49\linewidth]{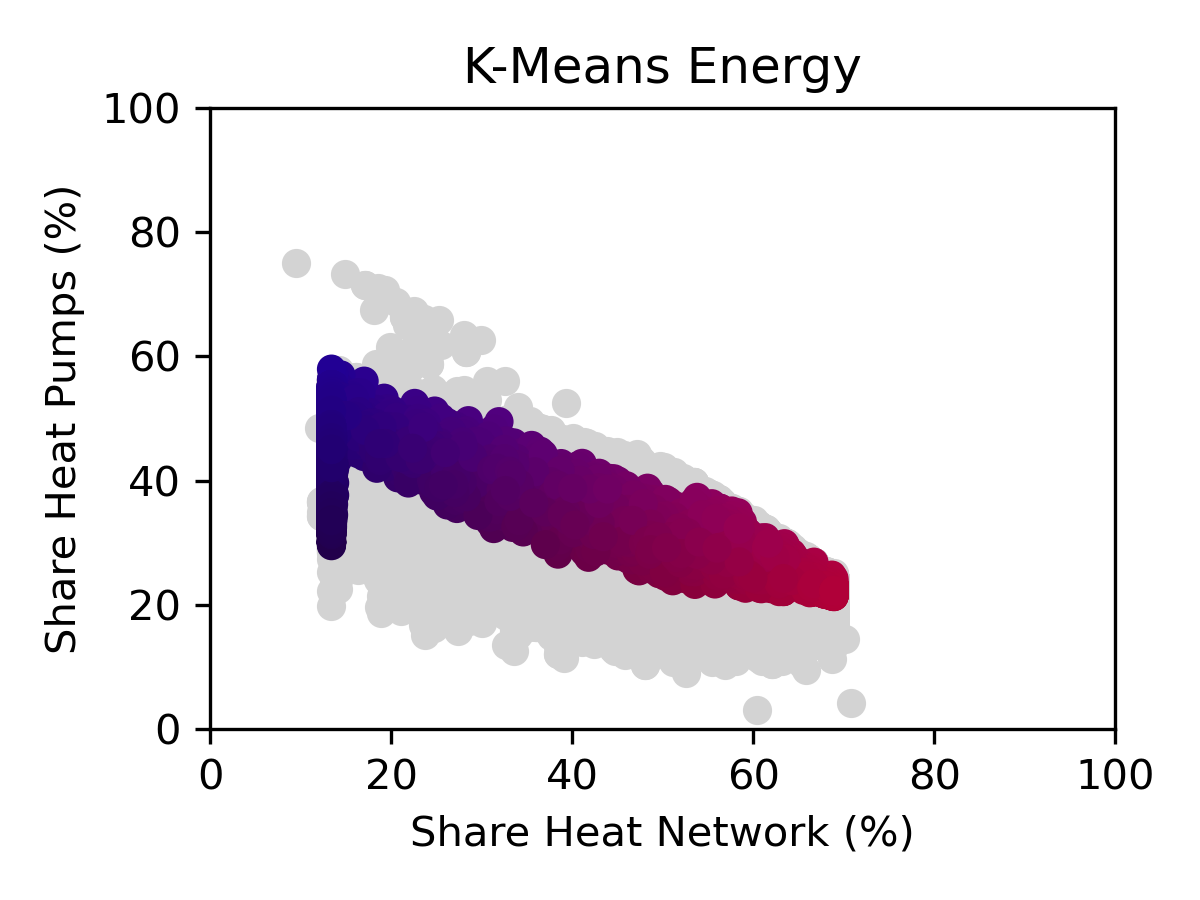}\\
    \includegraphics[width=0.49\linewidth]{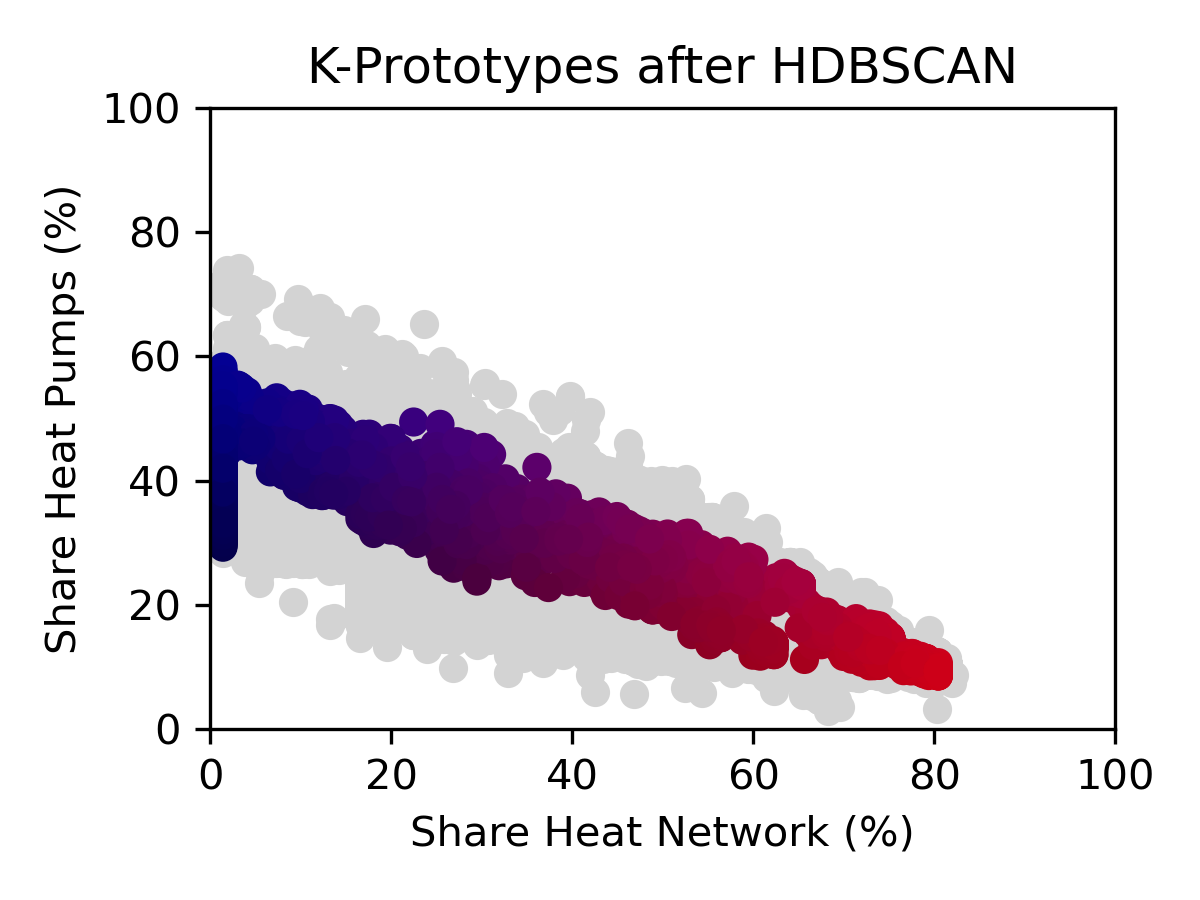}
    \includegraphics[width=0.49\linewidth]{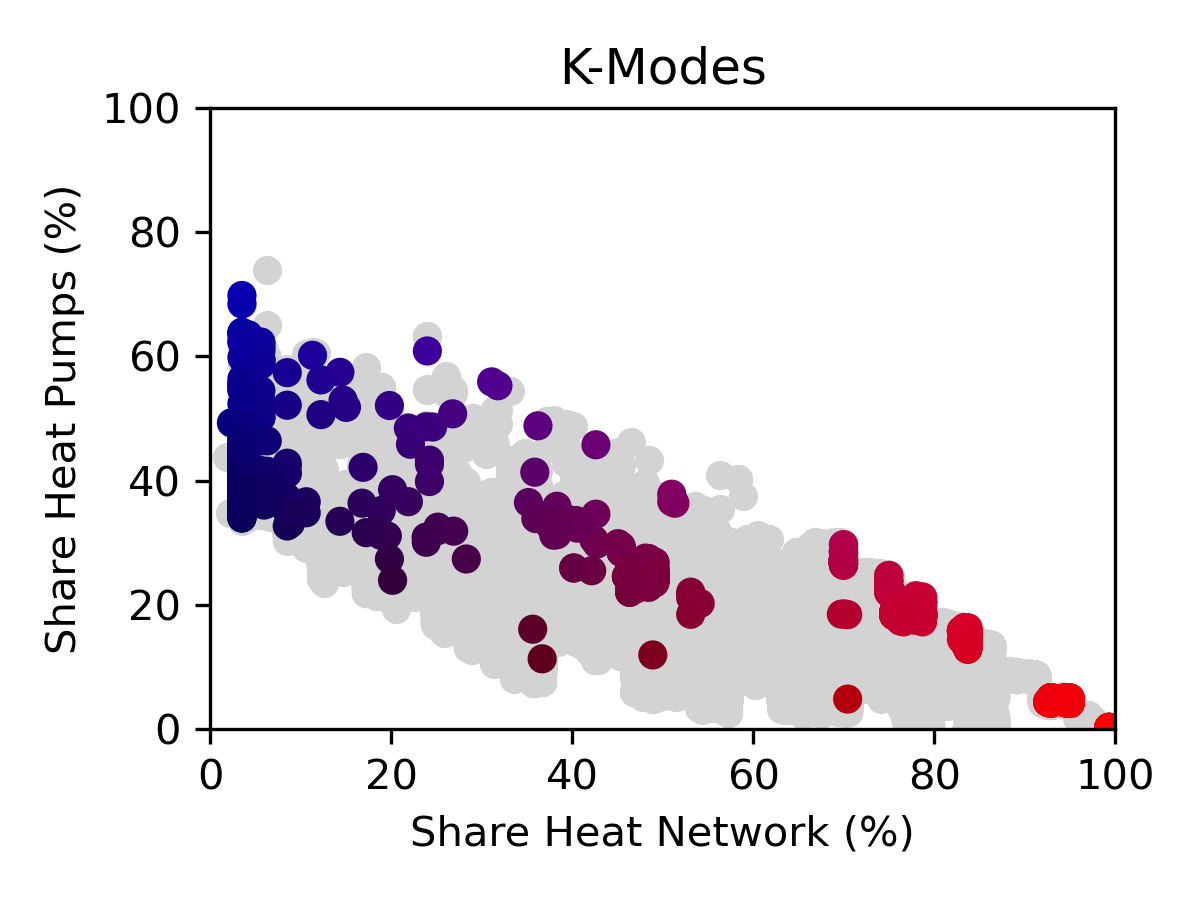}
    \caption{
        Sampled Pareto-plane in terms of global share of heating technologies
        for four different aggregation methods.
        Discarded solutions are displayed in gray,
        Pareto-optimal solutions are coloured.
        As the three add up to \SI{100}{\%}, the colour information
        (black: gas boiler, blue: heat pump, red: heat network)
        is redundant.
    }
    \label{fig: sampled optimal space}
\end{figure}

In the present study, heat can be provided by one of the three options,
gas boiler, heat pump, and heat network.
Fig.~\ref{fig: sampled optimal space} shows which area of the possible
combinations has been sampled by the Pareto-search.
It is clearly visible that one particular Pareto-optimal
solution has not been sampled by any of the runs:
In our model, not replacing any gas boiler and not installing any PV
would have resulted in no investment costs which is trivially Pareto-optimal.
\begin{figure}[thb]
    \centering
    \includegraphics[width=0.49\linewidth]{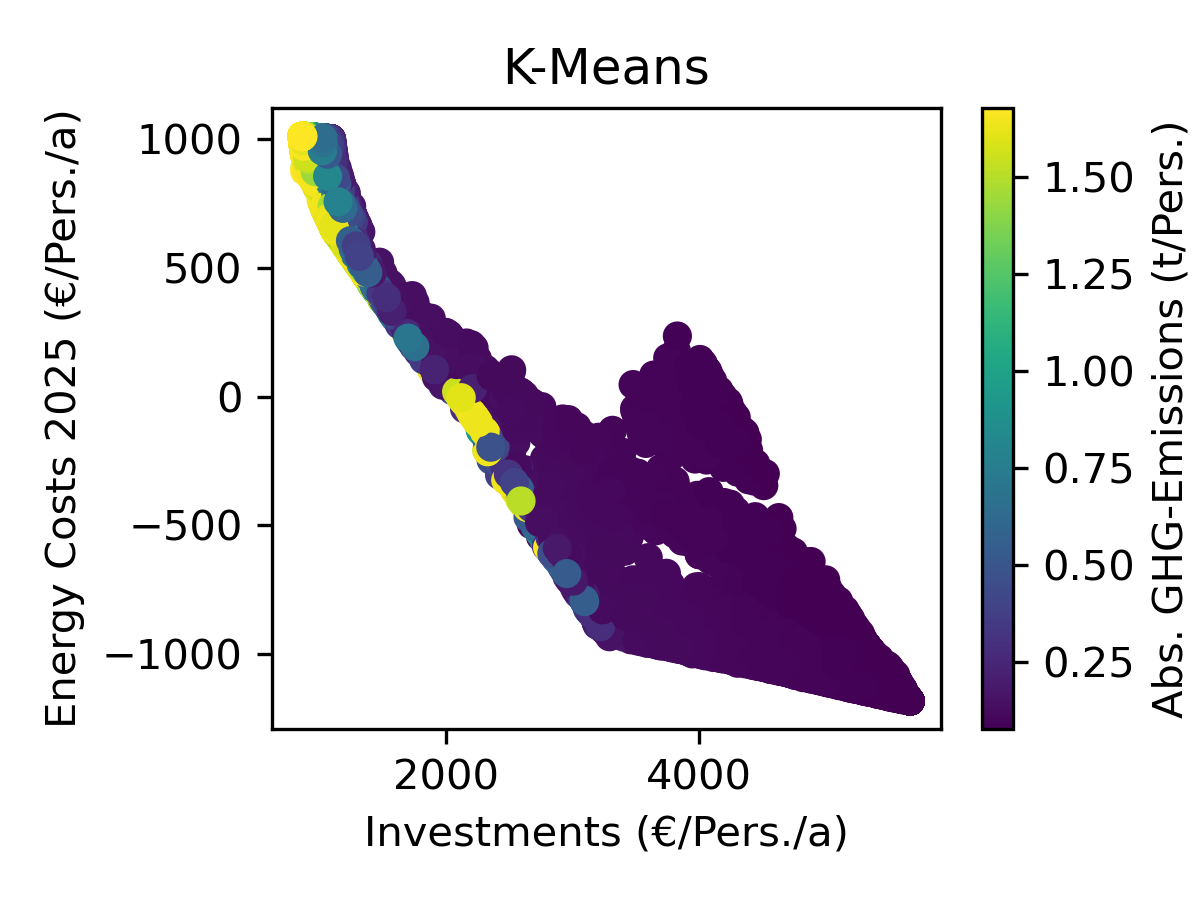}
    \includegraphics[width=0.49\linewidth]{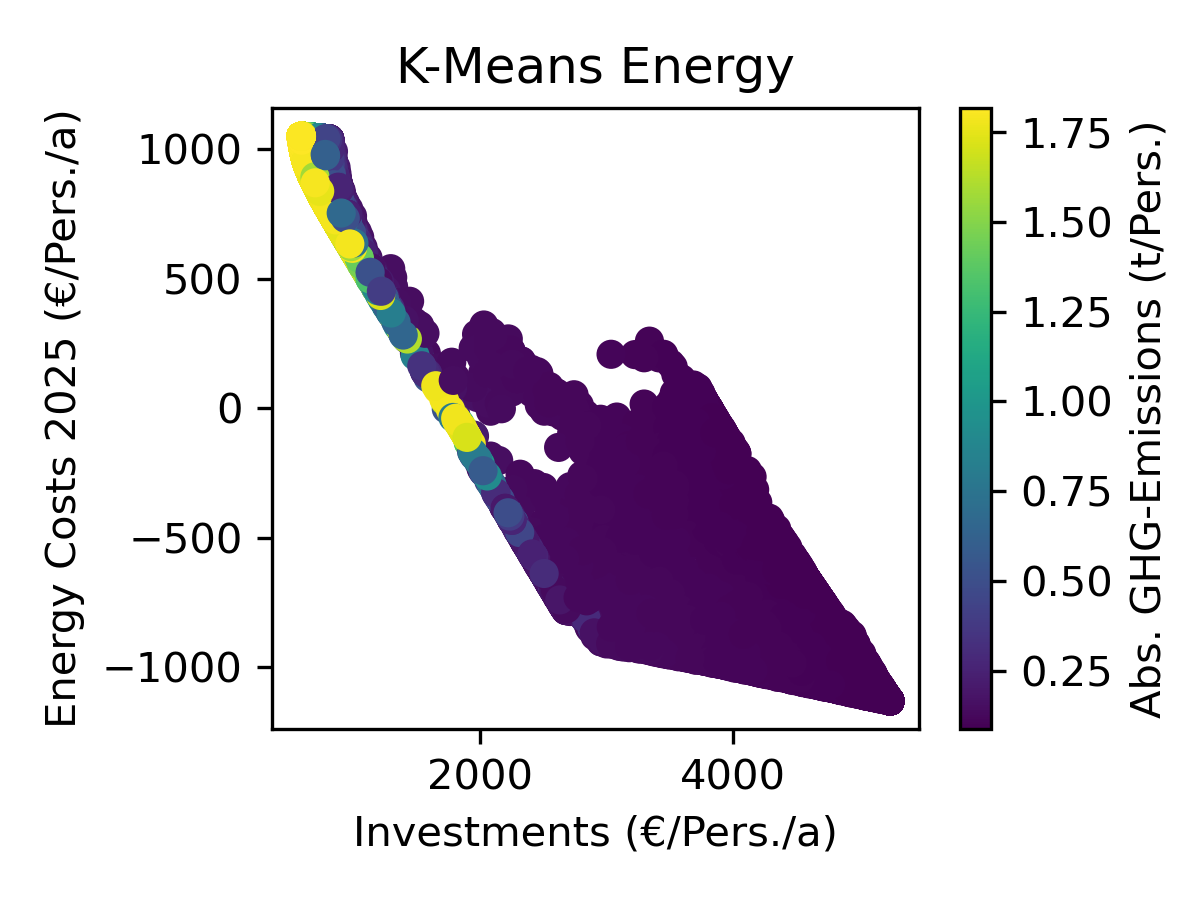}\\
    \includegraphics[width=0.49\linewidth]{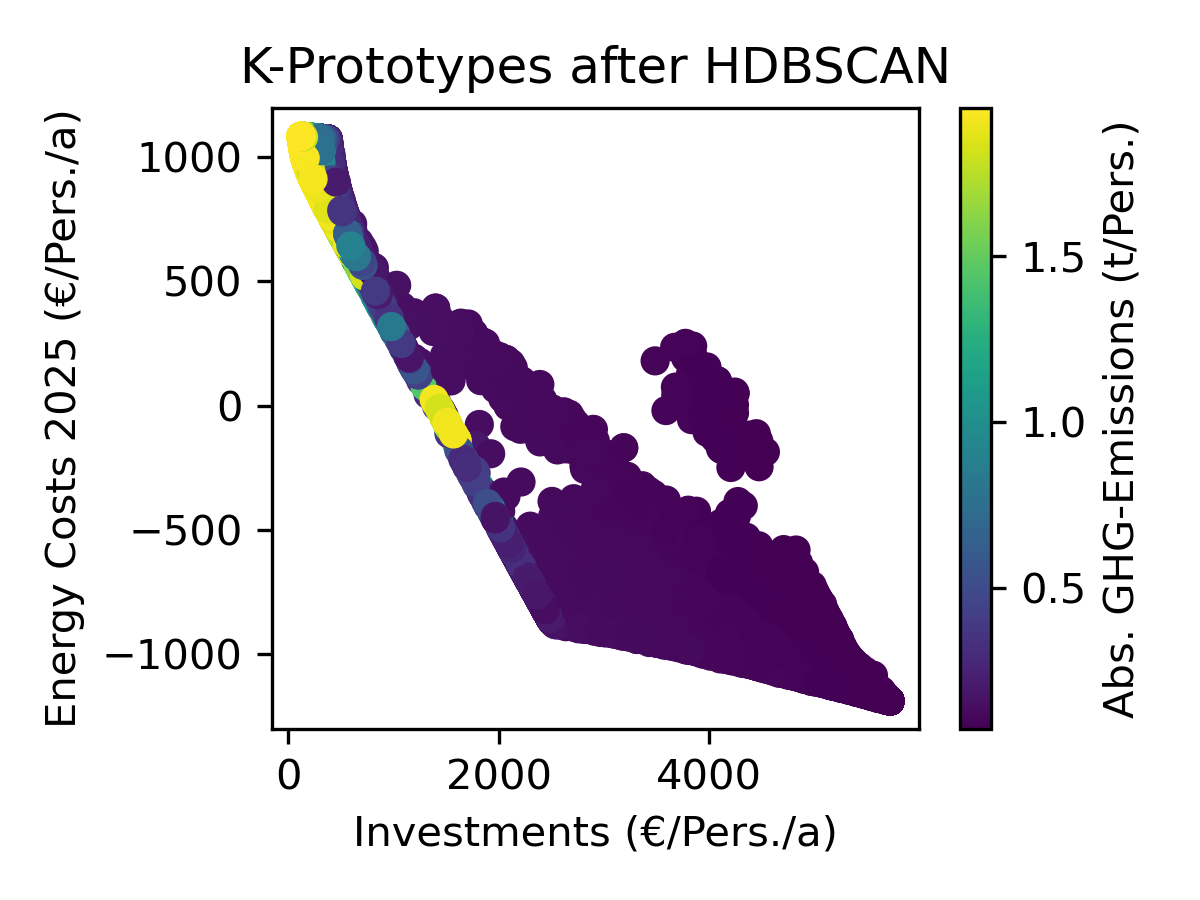}
    \includegraphics[width=0.49\linewidth]{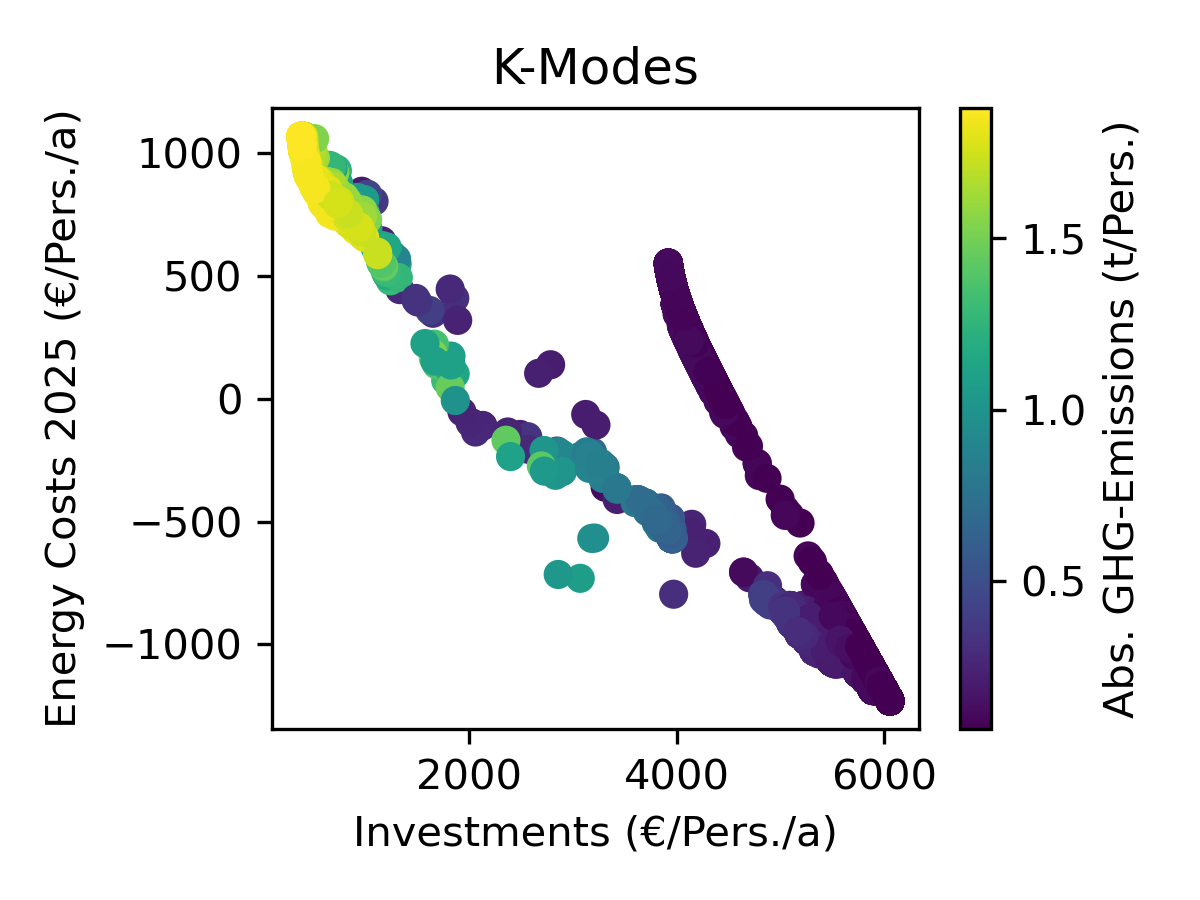}
    \caption{
        Pareto-plane in the space of conflicting goals.
        Starting from the low-investment edge,
        investment can first reduce both,
        operational costs and emissions.
        However, assuming energy costs of 2025,
        investments need to be subsidised to be profitable.
    }
    \label{fig: pareto optimal kpi}
\end{figure}
However,
when looking at the (in this case also three-dimensional) space of conflicting goals,
as displayed in Fig.~\ref{fig: pareto optimal kpi},
it can be seen that already the found solutions approach
rather low investment values.

In the comparison of the four clustering methods,
three produce similar aggregated configurations (Fig.~\ref{fig: sampled optimal space})
as well as similarly shaped Pareto-planes (Fig.~\ref{fig: pareto optimal kpi}).
From this perspective,
K-Modes deviates most, as it has a lower density of sampled points.
It does also not find the elbow at at the front of only investment
and operational costs.
At this one-dimensional front, the other methods show a point at about
\SI{-1000}{\sieuro/Pers./a} of energy costs,
where additional investment starts to reduce the energy costs
with a reduced slope.
One reason for this might be that configurations are made impossible
by the reduction of variables (cf. Tab.~\ref{tab: 10-5 variables}):
When the geographical groups strongly align with the presence of specific representative buildings,
it is there is decreased possibility to connect a heat network just in direct neighbourhoods.
The reason is that neighbouring buildings a different representative building type are often in a different geographical group
(cf. Fig~\ref{fig: map clusters}).
For purely geographical K-Means on the other extreme side, neighbouring buildings of the same representative building type are more likely to
be in different geographical groups compared to the methods that consider energy.
Thus, it is expectable that these neighbouring buildings
that would optimally take the same decision will be
treated differently.

Selecting concrete solutions that are comparable between the four clustering methods is not trivial.
In particular, the Pareto-front will always be sampled with a limited resolution.
This includes that particular weightings of optimisation goals,
especially extreme weightings like the cheapest solution,
might not be sampled -- as mentioned before.

\begin{figure}[thb]
    \centering
    \includegraphics[width=0.49\linewidth]{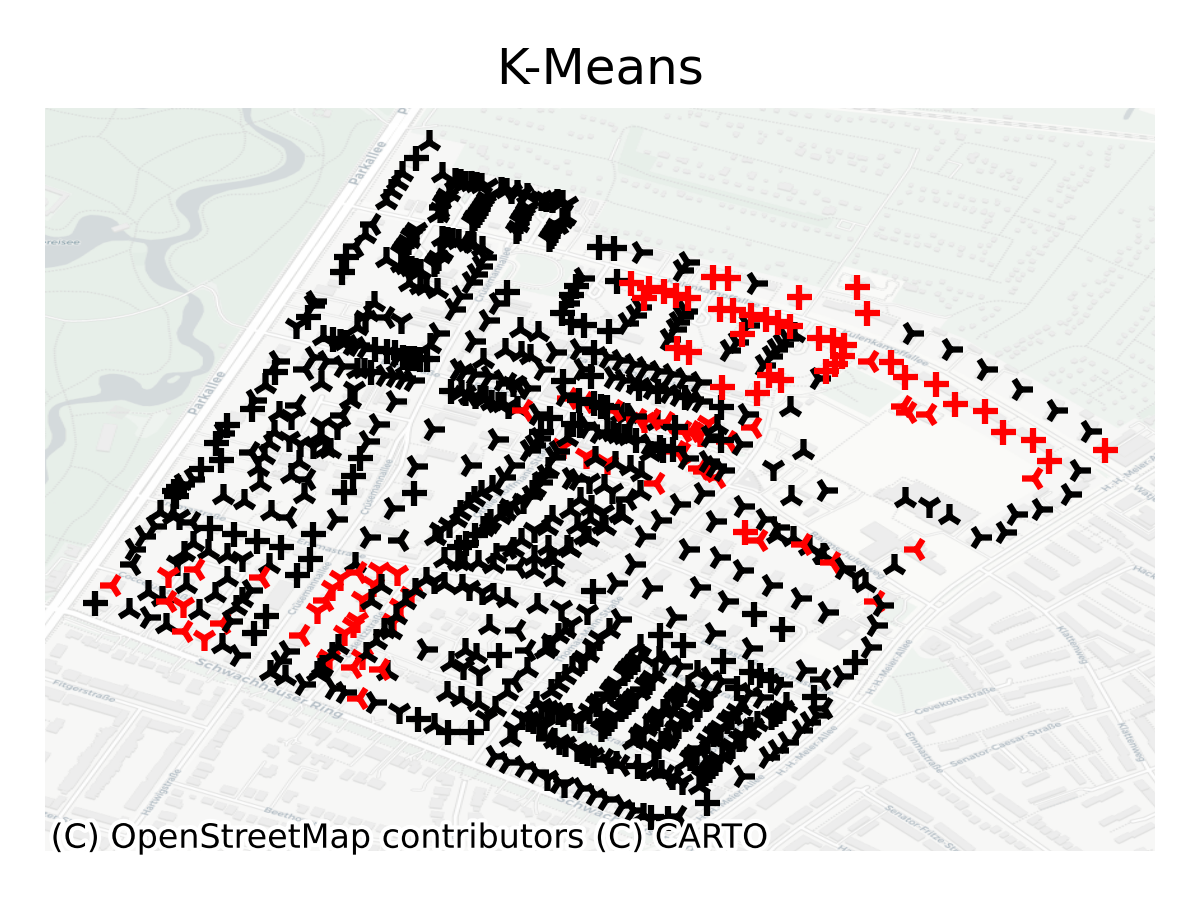}
    \includegraphics[width=0.49\linewidth]{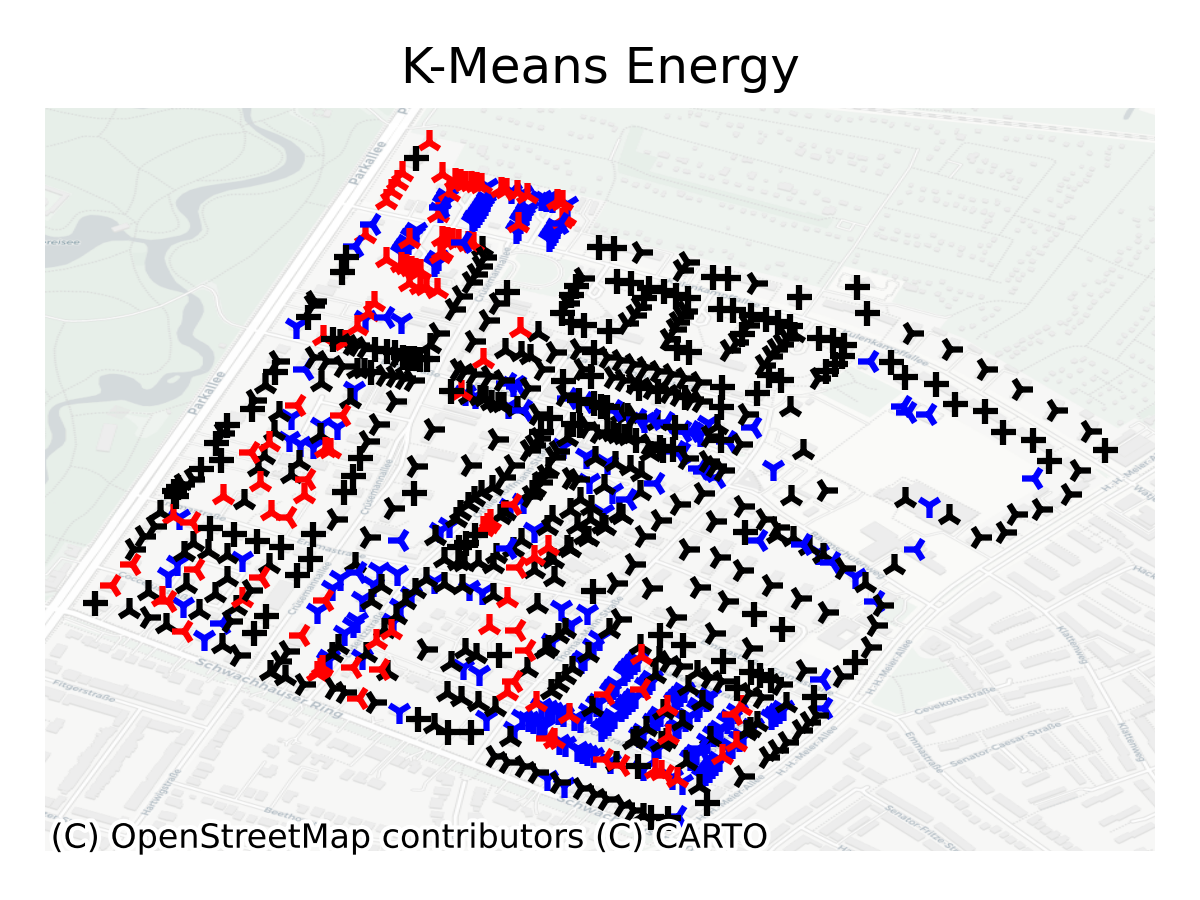}\\
    \includegraphics[width=0.49\linewidth]{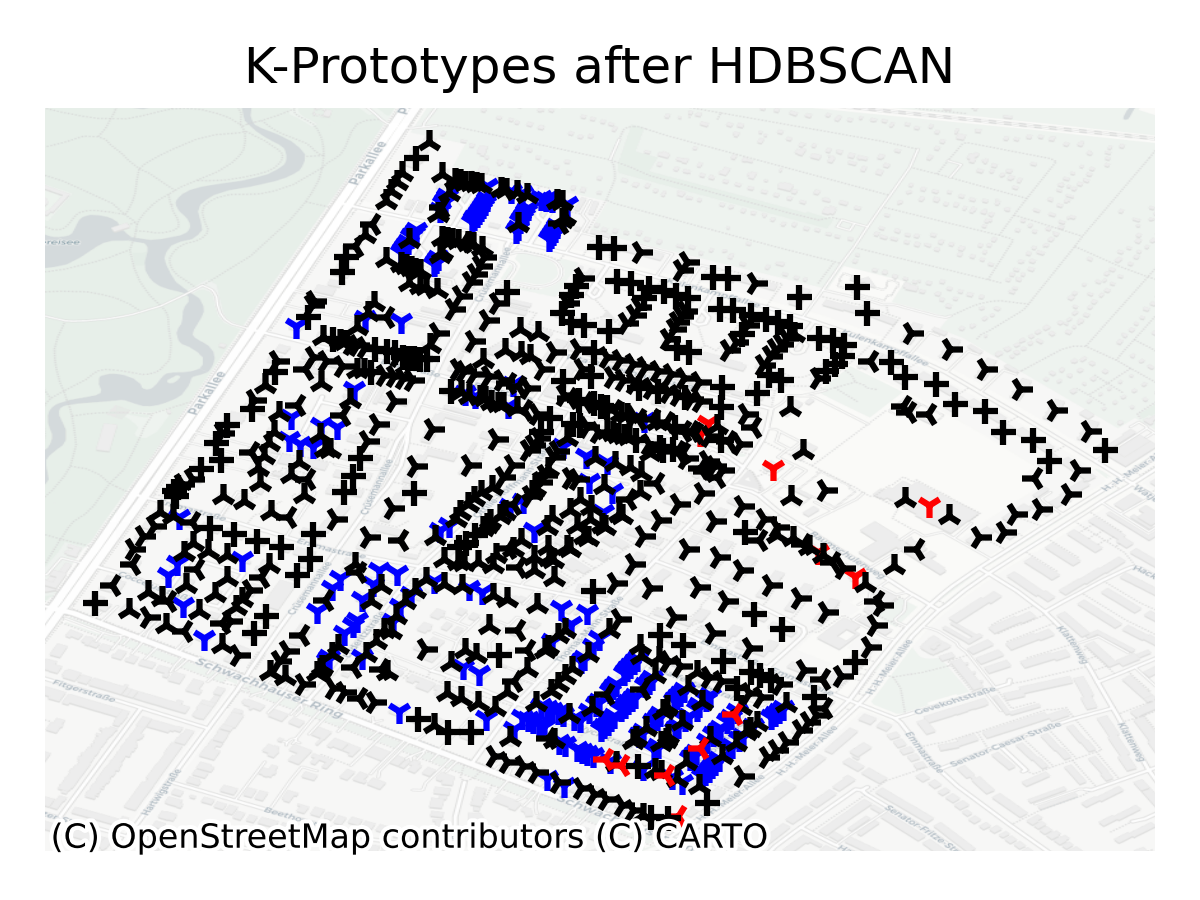}
    \includegraphics[width=0.49\linewidth]{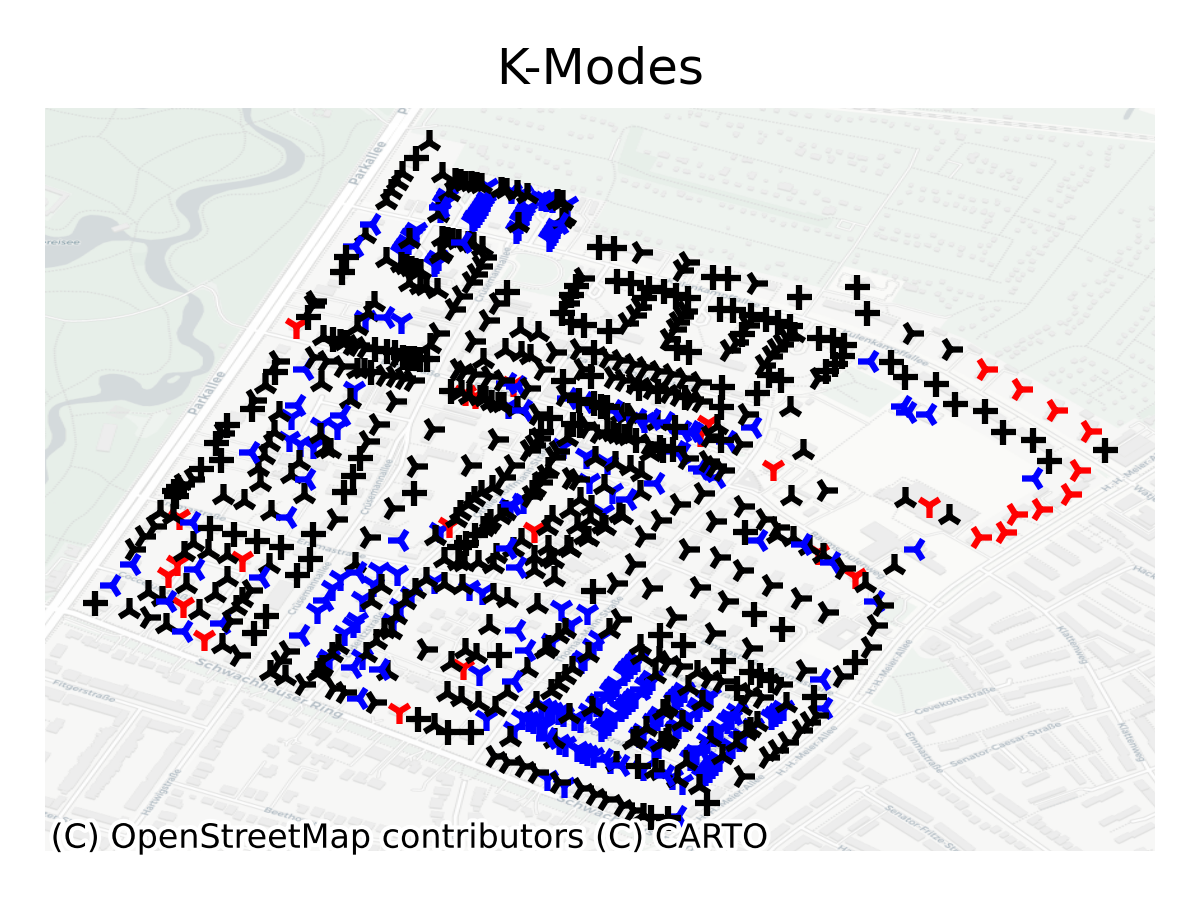}
    \caption{
        Technology choice in selected solutions based on grouping method (see Tab.\ref{tab: 500EUR-configs}).
        Heat network (red), heat pump (blue), gas boiler (black).
    }
    \label{fig: heat tech maps 500}
\end{figure}

First we select the solution with the lowest possible investment that
still has energy costs in 2025 that are below \SI{500}{\sieuro/Pers}.
In Fig.~\ref{fig: pareto optimal kpi},
this region looks densely sampled with a similarly shaped Pareto plane.
The filtered solutions are listed in Tab.~\ref{tab: 500EUR-configs}.
By design, all four configurations match the desired energy costs.
The lowest possible investment, however, deviates in a range between
\SI{756}{\sieuro/Pers./a} and \SI{1339}{\sieuro/Pers./a}.
For the emissions, two of the solutions show values
of approximately \SI{1.2}{t/Pers.} and two of approximately \SI{1.7}{t/Pers.}.
In total, the solution found by K-Means Energy can be considered the best,
K-Prototypes is better than K-Means (without energy information),
but there is no strictly worst solution.
As we searched for low investment costs,
gas boilers have the biggest share of all heating technologies.
The share goes from \SI{87.8}{\%} (K-Means) down to \SI{57.6}{\%} (K-Modes).
All but K-Means (not heat-pumps at all) prefer heat-pump over heat-network.
The choice per building is visualised in Fig.~\ref{fig: heat tech maps 500}.
\begin{table}
    \centering
    \begin{tabular}{c c|c c c c}
        Quantity & Unit & K-Means & K-Means & K-Prototypes & K-Modes \\
         &  & & Energy  & after HDBSCAN &  \\
        \hline
        Heat Network & \si{bldg}      & 105  & 115 &  12 &  30\\
        Heat Pump & \si{bldg}         &   0  & 249 & 197 & 258\\
        Gas Boiler & \si{bldg}        & 754  & 495 & 650 & 571\\
        Investments & \si{\sieuro/Pers/a} & 1339  & \textbf{756}  & 1211	& 1133\\
        Energy Costs 2025 & \si{\sieuro/Pers}
                                    & 499	& 497  & 499	& 500\\
        Emissions  & \si{t/Pers}    & 1.64	& \textbf{1.17}	& \textbf{1.16}	& 1.70\\
    \end{tabular}
    \caption{
        Configurations with lowest investment costs that have a limit
        of \SI{500}{\sieuro/Pers.} of energy costs in 2025.
        Best and close to best values are bold.
    }
    \label{tab: 500EUR-configs}
\end{table}
The three of the algorithms that place decentralised heat pumps,
show a high overlap (183 buildings assigned consistently).
Assignment of heat network, however, is very inconsistent with only
two buildings being attached to the heat network in three out of four cases.

\begin{figure}[thb]
    \centering
    \includegraphics[width=0.49\linewidth]{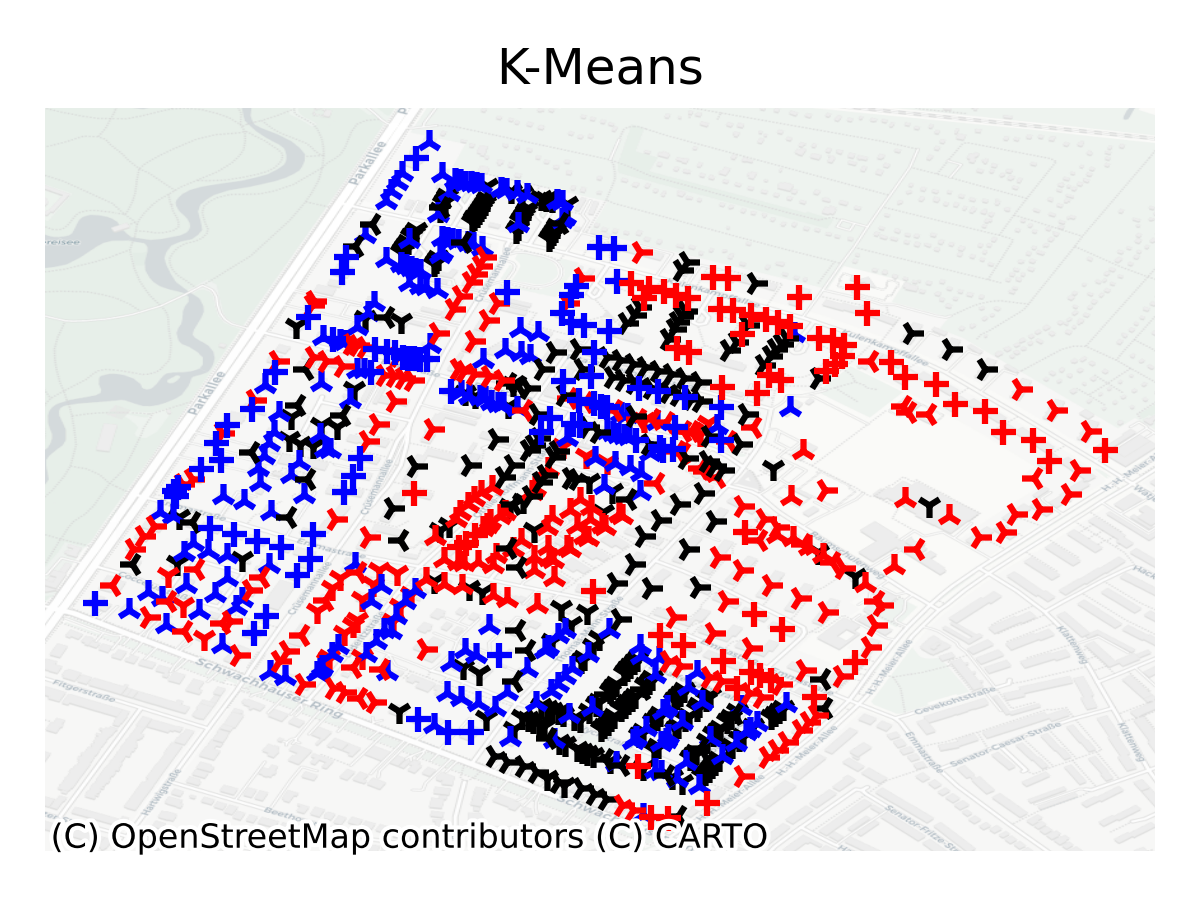}
    \includegraphics[width=0.49\linewidth]{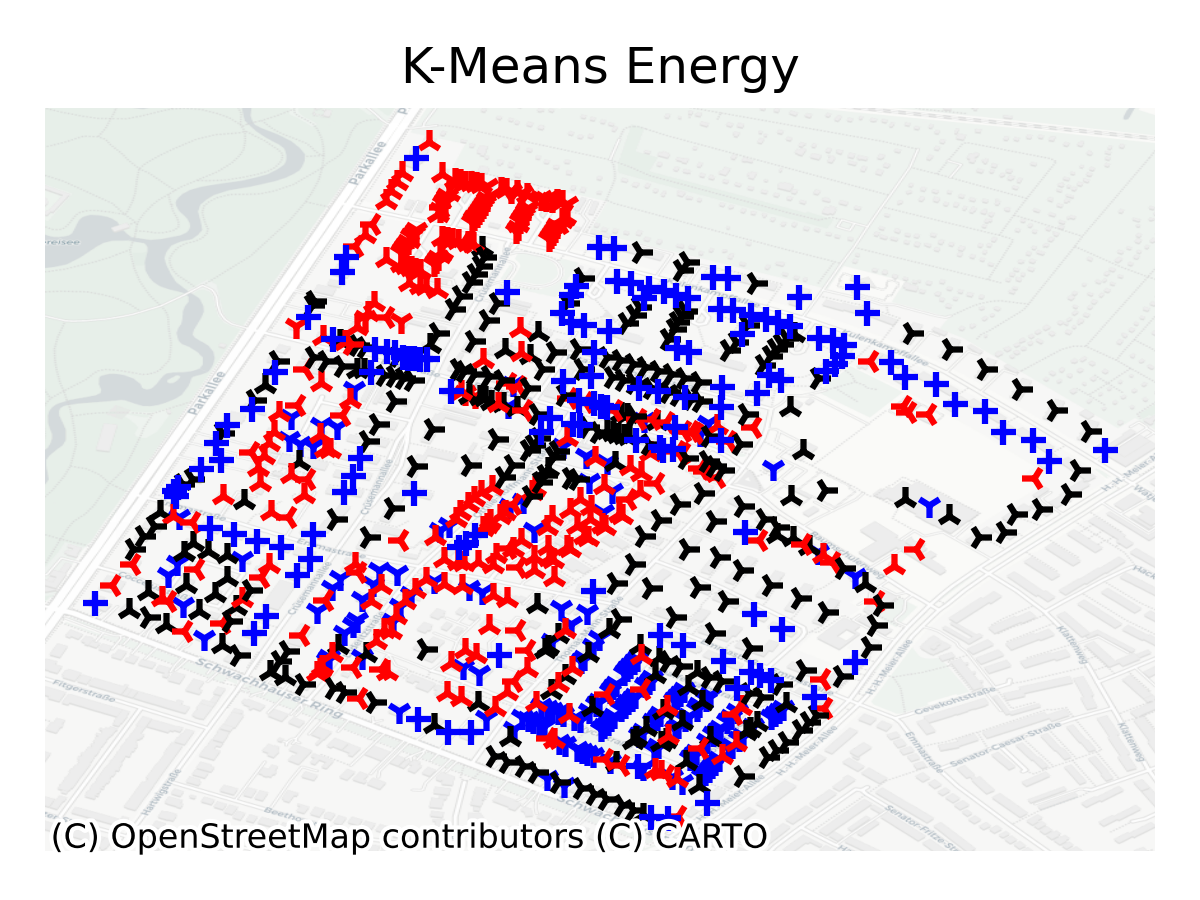}\\
    \includegraphics[width=0.49\linewidth]{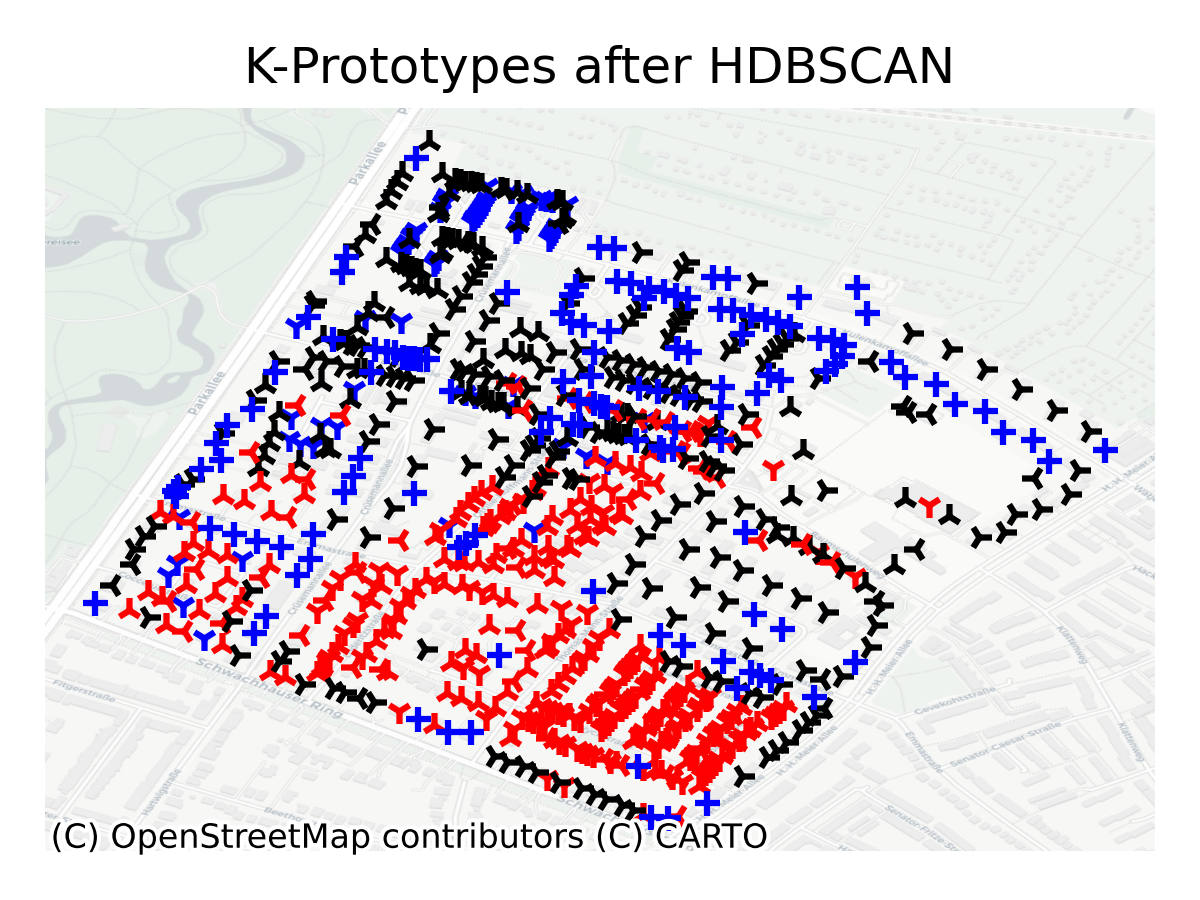}
    \includegraphics[width=0.49\linewidth]{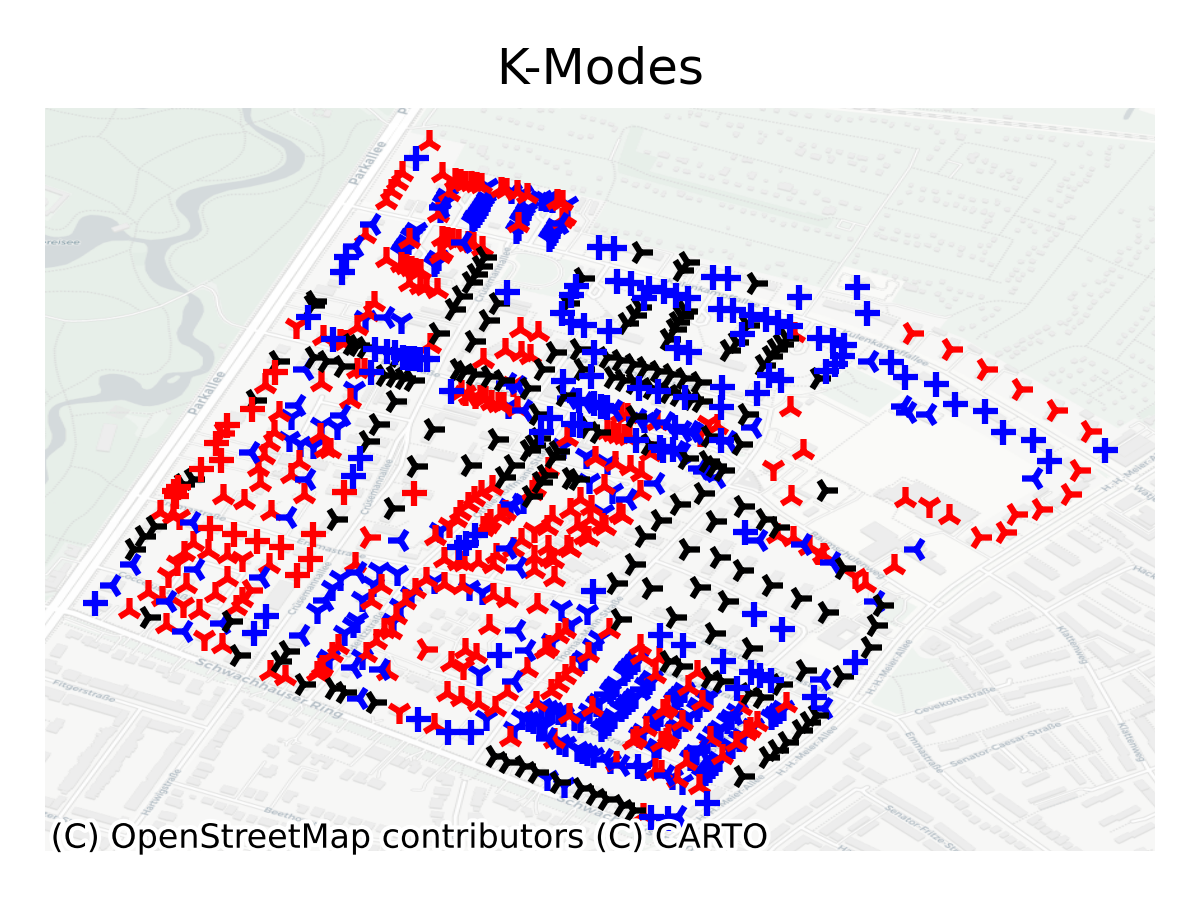}
    \caption{
        Technology choice in selected solution based on grouping method (see Tab.\ref{tab: 33-configs}).
        Heat network (red), heat pump (blue), gas boiler (black).
    }
    \label{fig: heat tech maps 33}
\end{figure}

Second, we select the configuration with a minimum of \SI{33}{\%} of buildings using both,
heat network and heat pump,
that has the lowest investment annuity.
We chose this specific configuration to find out which areas and buildings
prefer which technology in a more or less equally distributed system.
The filtered results are summarized in Tab.~\ref{tab: 33-configs}.
It can be seen that the configurations show similar investment costs
(\SI{2572}{\sieuro/Pers./a} to \SI{2924}{\sieuro/Pers./a})
and emissions (\SI{0.11}{t/Pers.} to \SI{0.22}{t/Pers.}),
but deviate in the energy costs in 2025 with values ranging
from a revenue of \SI{339}{\sieuro/Pers.} to spendings of \SI{138}{\sieuro/Pers.}.
Interestingly, there is a tendency that the solutions with slightly lower investment costs
tend to prefer lower penetration with gas boilers.
Still, the filter makes sure that
all three technologies are used in a significant number of buildings,
as also seen in Fig.~\ref{fig: heat tech maps 33}.

\begin{table}
    \centering
    \begin{tabular}{c c|c c c c}
        Quantity & Unit & K-Means & K-Means & K-Prototypes & K-Modes \\
         &  & & Energy  & after HDBSCAN &  \\
        \hline
        Heat Network & bldg      & 285 & 290 & 348 & 311 \\
        Heat Pump & bldg         & 269 & 289 & 199 & 372 \\
        Gas Boiler & bldg        & 305 & 280 & 312 & 176 \\
        Investments & \si{\sieuro/Pers/a} & 2924 & \textbf{2642} & \textbf{2572} & 2783 \\
        Energy Costs 2025 & \si{\sieuro/Pers}
                                    &  -86 & \textbf{-339} & \textbf{-412} & 138 \\
        Emissions  & \si{t/Pers}    & \textbf{0.11} & \textbf{0.12} & \textbf{0.12} & 0.22 \\
    \end{tabular}
    \caption{
        Configurations with lowest investment that have at least \SI{33}{\%}
        of both, heat network and heat pump.
        Best and close to best values are bold.
    }
    \label{tab: 33-configs}
\end{table}

When comparing the three maps,
it can be seen that some of the buildings show a clear preference over
all four methods, while others are assigned different technologies between the runs.
This is independent from the assigned archetype building.
A numeric evaluation yields that
\SI{10.4}{\%} of the buildings are consistently assigned a gas boiler,
\SI{6.4}{\%} a heat pump, and
\SI{5.6}{\%} a heat network connection
across all three of the favoured methods.
While it sounds underwhelming that only less then a third
of the buildings have the heating technology chosen independent
from the grouping method,
the numbers are still significantly bigger compared to
purely random assignment
(\SI{5.0}{\%}, \SI{3.4}{\%}, and \SI{2.4}{\%}, respectively).

The sensitivity of the result on the used clustering method can be interpreted in two different ways:
It is possible that the clustering can introduce a considerable error,
but it is also possible that the selected region for the study is so homogenous
that even small inaccuracies change the result significantly.
In the latter case, the deviation between the results would imply
that the choices do not make a big difference.
The overall shape of the Pareto-space as depicted in Fig.~\ref{fig: pareto optimal kpi}
supports this interpretation at least for three of the investigated methods.

\section{Conclusion}

This study presented a method to preprocess and aggregate geodata
to be able to optimise the energy system configuration of an ensemble
of a large number of buildings.
After the selection of the area,
all presented steps can in be applied in an automatic workflow without any manual input.
While details of the actual results depend on simulation parameters,
the overall performance of investigated clustering methods is comparable.
In our study, methods that balance influence of geographic location
and energy properties (K-Means Energy and K-Prototypes after HDBSCAN)
showed the best performance.
In particular, not considering energy at all when building geographical groups
hides options for very small heating grids just between similar buildings
in close proximity.
While the low number of variables was a disadvantage here,
it should be noted that K-Modes might still be advantageous
if a more complex energy system model was used.
Especially when including more representative buildings,
storage and optimised energy-exchange between the buildings,
the overall runtime might increase to an extent,
that the expected advantage in run-time might outweigh
the reduced number of options.

We made the choice to use building parameters normalised by the
building footprint to allow for a low number of representative buildings
all representing a similar number of real buildings.
In particular, we also normalised the costs of a connection to a heating grid.
The alternative choice of non-normalised buildings
would allow for non-linear investment costs.
If the model still stays in a linearisable regime,
even the argument about commutable operation from the introduction stays applicable.
However, better representation with lower number of representatives seems
to be more important at the moment.
Note that this argument might be invalidated
if other highly non-linear investment options come into scope,
resulting in a higher number of variables for the normalised representation.

For future work,
we suggest to test the method including the different grouping methods with another region
to assess whether an area with clear advantages of one technology is correctly identified.
This might be at the edge of a town where a densely populated area transitions
into a more rural area.
Also, the method to use representative roofs should be validated and possibly fine-tuned.
Lastly, it would be interesting to include additional
demands and technologies, esp. cooling demand and bidirectional heat networks.

\section*{Declarations}

\subsection*{Acknowledgements}
The authors thank ``wesernetz'' for the provision of the real-world
heat demand data used in the use case scenario.
They also acknowledge Diana Maldonado for contributing
code for the Pareto-optimisation and for valuable discussions,
in particular about the selection of representative configurations.

\subsection*{Author contributions}
Patrik Schönfeldt: conceptualization, methodology, software, formal analysis, investigation, writing -- original draft;
Elif Turhan: methodology, writing -- reviewing and editing.

\subsection*{Funding}
This work has been funded by the German Federal Ministry of Research, Technology and Space as part of the Project ``Wärmewende Nordwest'', Grant No. 03SF0624L.

\subsection*{Data availability}
The current study uses heat data
received under a no-disclosure agreement.
Conditions apply to access the data.

\begin{appendices}

\section{Two-step clustering scan}
\label{sec: two-step scan}

As described in Sec.~\ref{sec: clustering},
we implemented a two-step clustering method,
where a data-driven definition of representative
buildings is followed by a geographical grouping of these.
Details on the applied grouping methods
can be found in~\ref{sec: geographical grouping}.

\begin{figure}[pthb]
    \centering
   \includegraphics[width=0.49\linewidth]{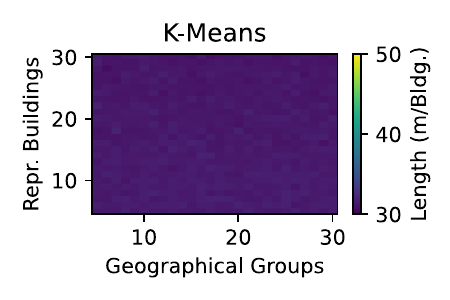}
   \includegraphics[width=0.49\linewidth]{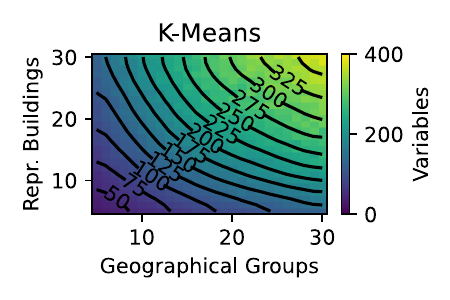}\\
   \includegraphics[width=0.49\linewidth]{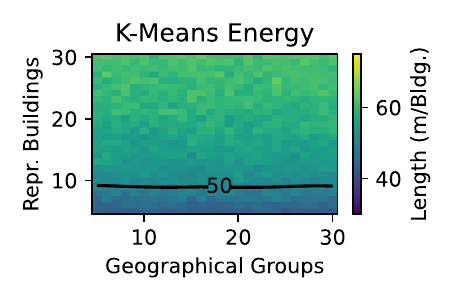}
   \includegraphics[width=0.49\linewidth]{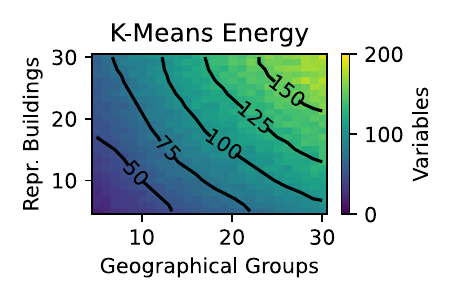}\\
   \includegraphics[width=0.49\linewidth]{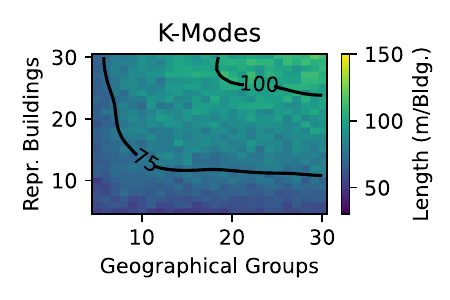}
   \includegraphics[width=0.49\linewidth]{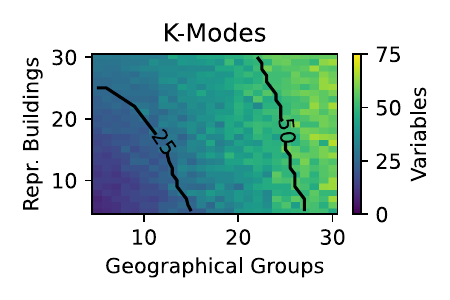}\\
   \includegraphics[width=0.49\linewidth]{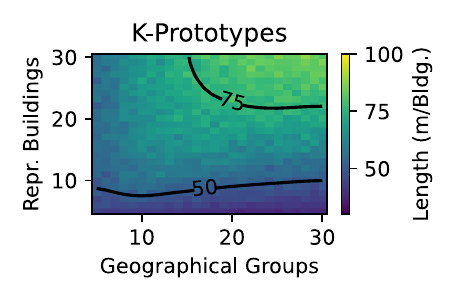}
   \includegraphics[width=0.49\linewidth]{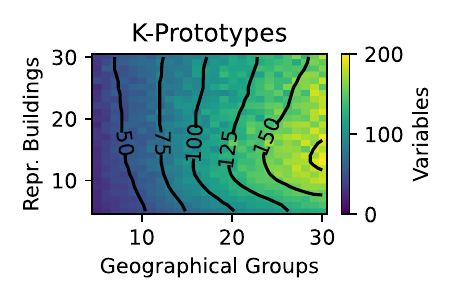}\\
   \includegraphics[width=0.49\linewidth]{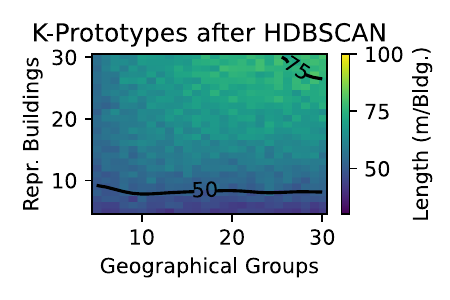}
   \includegraphics[width=0.49\linewidth]{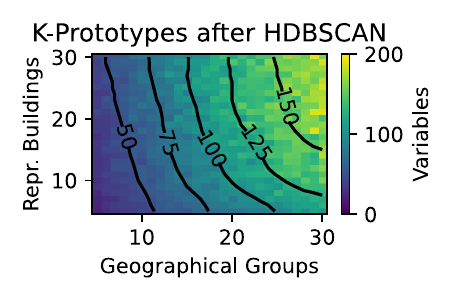}
    \caption{
        Average single-group line length (left)
        and number of variables (right).
        We introduced smoothened lines to guide the eye.
        Note that these are not \SI{100}{\%} accurate.
    }
    \label{fig: performance of clustering algorithms}
\end{figure}

Here, we present the result of a scan over the number
of both, representative buildings and geographical groups,
show how these effect the number of variables needed to
describe the final model and the average single-group line length,
which we take as a measure for geographical cohesion of the grouping method.
The results are visualised in
Fig.~\ref{fig: performance of clustering algorithms}.
In the figure,
the right hand side shows the average length of a heat network,
if heat networks were built connecting all buildings in the same geographical group.
The left hand side displays the number of existing combinations
between building category and geographical group.
This translates to the number of variables needed in an optimisation model.

Depending on the method used for geographical grouping,
increasing the number of these groups does not decrease the
geographical cohesion.
For both K-Prototypes approaches and K-Modes this is
no longer true if the number of representative buildings
is increased too much for the targeted number of geographical groups.
The absolute number of variables is clearly highest
for purely geographical K-Means.
This is no surprise as it has no mechanism
to group similar buildings.
K-Means including energy as well as both approaches using K-Prototypes show similar performance in this indicator,
with an advantage of K-Prototypes when the number
of representative buildings is increased for a constant number of geographical groups.
By favourably grouping the same representative buildings in one geographical group,
K-Modes manages to keep the number of variables lowest.
For K-Prototypes without HDBSCAN,
we observe that increasing the number of representative buildings
for high number of geographical groups can decrease the number of variables.
As the total number of variables is still relatively high,
we did not investigate this effect further.

\section{Energy system model}
\label{sec: energy system model}

For the present study, we apply an energy balance model without any storage.
The time series for heat demand, PV potential, and the COP curve of the air source
heat pump are calculated from the typical meteorological year for the area (TMY of 2015, tile 530876088493) as published by the DWD (Deutscher Wetterdienst, German Weather Service)~\cite{DWD_try2015}.
Using this weather time series,
annual energy demand is converted to energy time series using
the software \texttt{demandlib}~\cite{demandlib_v022}.
We also use the air temperature of the TMY
as the source temperature \(T_\mathrm{S}(t)\)
to calculate the efficiency of the decentralised air source heat pumps as
\begin{equation}
            1/\eta_\mathrm{C}(t) = c_\mathrm{pf}
                \times \frac{T_\mathrm{F}}{T_\mathrm{F} - T_\mathrm{S}(t)}
\end{equation}
assuming a constant thermodynamic efficiency of \(c_\mathrm{pf} = 0.45\),
and a global flow temperature of \(T_\mathrm{F} = \SI{328.15}{K}\).
The irradiation data can be used to calculate the solar power potential time series
using PVLib~\cite{Anderson2023} and Eq.~\eqref{eq: representative roof weights}
as
\begin{equation}
    P_{\mathrm{PV},b}(t) = \frac{A_b}{\SI{1}{m^2}} \sum_{n=0}^{N-1} \omega_{n,b} \times P_{n}(t).
\end{equation}

The building energy balance depends on the decision of the heating technology.
We define
\(\delta_{\mathrm{HP},b,g} = 1\)
if representative building \(b\) in group \(g\) uses a heat pump,
\(\delta_{\mathrm{GB},b,g} = 1\) for a present gas boiler,
and \(\delta_{\mathrm{HN},b,g} = 1\) if the building is connected to the heat network.
In our present model, the heat supply options are mutually exclusive, meaning
\begin{equation}
    \delta_{\mathrm{HP},b,g} + \delta_{\mathrm{GB},b,g} + \delta_{\mathrm{HN},b,g} = 1.
\end{equation}
PV investment, on the other hand, is assumed to be continuous \(0 \le q_{\mathrm{PV},b,g} \le 1\).
From this, the time series for the building energy balance can be calculated as
\begin{subequations}
\begin{equation}
    P_{\mathrm{el},b,g}(t) = P_{\mathrm{ESLP},b}(t) + \frac{\delta_{\mathrm{HP},b,g}}{\eta_{\mathrm{C},b}(t)} P_{\mathrm{HSLP},b}(t) - q_{\mathrm{PV},b,g} P_{\mathrm{PV},b}(t),
\end{equation}
for electricity and
\begin{equation}
    P_{\mathrm{gas},b,g}(t) = \delta_{\mathrm{GB},b,g}\times P_{\mathrm{HSLP},b}(t)
\end{equation}
for gas.
The remaining power time series for the central heat pump is
\begin{equation}
    P_{\mathrm{el,HN}}(t) = \frac{1}{\eta_{\mathrm{C},b}(t)} \sum_{b\in\mathfrak{B}}\sum_{g \in \mathfrak{G}} \delta_{\mathrm{HN},b,g} P_{\mathrm{HSLP},b}.
    \label{eq: power central heat pump}
\end{equation}
\end{subequations}
Further, for the emissions as defined in Eq.~\eqref{eq: indicator emissions},
the global energy balance is relevant.
We have
\begin{subequations}
    \begin{equation}
        P_{\mathrm{el}}(t) = \sum_{b\in\mathfrak{B}}\sum_{g \in \mathfrak{G}} P_{\mathrm{el},b,g}(t) + P_{\mathrm{el,HN}}(t)
    \end{equation}
for electricity, and
    \begin{equation}
        P_{\mathrm{gas}}(t) = \sum_{b\in\mathfrak{B}}\sum_{g \in \mathfrak{G}} P_{\mathrm{gas},b,g}(t)
    \end{equation}
for gas.
\end{subequations}

Using the corresponding power time series,
the energy costs in 2025 of representative building \(b\) in group \(g\)
as used in Eq.~\eqref{eq: indicator energy costs 2025} is
\begin{subequations}
\begin{align}
    C_{b,g} =
        &c_\mathrm{el}\int \max\left(
            P_{\mathrm{el},b,g}(t), \SI{0}{kW}
        \right)\;\mathrm{d}t\nonumber\\
        - &c_\mathrm{feedin} \left|
            \int \min\left(P_{\mathrm{el},b,g}(t), \SI{0}{kW}\right) \;\mathrm{d}t
        \right|\nonumber\\
        +  &c_\mathrm{gas} \int P_{\mathrm{gas},b,g}(t)\;\mathrm{d}t.
        \label{eq: building enery costs}
\end{align}
The respective costs are
\(c_\mathrm{el} = \SI{0.3}{\sieuro/kWh}\),
\(c_\mathrm{feedin} = \SI{0.08}{\sieuro/kWh}\), and
\(c_\mathrm{gas} = \SI{0.1}{\sieuro/kWh}\).
The energy costs for the heat network \(C_{\mathrm{el},HN}\) are assumed to be proportional to its energy demand
\begin{equation}
    C_{\mathrm{el},HN} = c_\mathrm{el} \int P_{\mathrm{el,HN}}(t)\;\mathrm{d}t.
\end{equation}
\end{subequations}
Assuming the same electricity costs for households
and the heat network operator is a simplification made
to make the results easier to understand focussing on the
differences due to the aggregation.
It also compensates for the optimistically high COP
when it comes to the cost calculation.

Considering the investment cost,
to allow for more flexibility when deciding about algorithms
for the heuristic Pareto-optimisation,
investment decision is not taken to be discrete.
Instead, we introduce
\(0 \le q_{\mathrm{HP},b,g} \le 1\),
and
\(0 \le q_{\mathrm{HN},b,g} \le 1\)
as continuous decision variables.
They are rounded to obtain the energy system model values
\begin{subequations}
    \begin{equation}
        \delta_{h,b,g} = \lfloor q_{h,b,g} + 0.5 \rfloor
    \end{equation}
for \(h \in \{\text{HP, HN}\}\).
For the actual investment, we then use the upper limit of the two
    \begin{equation}
        \hat{q}_{h,b,g} = \max(q_{h,b,g}, \delta_{h,b,g}),
    \end{equation}
so that there is an incentive for the heuristic to converge to
single-technology solutions.
\end{subequations}
Finally, to compute the investment as defined in Eq.~\eqref{eq: indicator investment costs},
we apply price scaling curves using values of~\cite{kww_technikkatalog}.
\begin{subequations}
    The costs of building heat pumps then are
    \begin{equation}
        C_{\mathrm{HP},b,g} =\hat{q}_{\mathrm{HP},b,g}\times \frac{\hat{P}_{\mathrm{heat}}^{0.705}}{\SI{1}{kW}}\times \frac{\SI{3830.5}{\sieuro}}{\SI{20}{a}}.
        \label{eq: HP costs}
    \end{equation}
    For peak powers \(\hat{P}_\mathrm{heat} \ge \SI{50}{kW}\),
    specifically for the heat network,
    the constants are first changed to \(0.793\) and \SI{3194.6}{\sieuro},
    and from \SI{300}{kW} on to \(0.755\) and \SI{1352.0}{\sieuro}.
    For PV systems we assume    
    \begin{equation}
        C_{\mathrm{PV},b,g} = P_{\mathrm{PV},b,g}\times \frac{\SI{1500}{\sieuro/kW}}{\SI{20}{a}}.
    \end{equation}
    The costs of the connection of a building to the heat network are
    \begin{equation}
        C_{\mathrm{HN},b,g} = \frac{\hat{q}_{\mathrm{HN},b,g}}{\SI{25}{a}}\times \left(\frac{\hat{P}_{\mathrm{heat}}}{\SI{1}{kW}}\times \SI{13.4}{\sieuro} + \SI{13976}{\sieuro} \right).
    \end{equation}
    Finally, the contribution of the heat network infrastructure itself is
    \begin{equation}
        C_{\mathrm{HN}} = l_\mathrm{HN} \times \frac{\hat{P}_{\mathrm{heat}}^{0.2029}}{\SI{1}{kW}}\times \frac{\SI{432.69}{\sieuro}}{\SI{50}{a}} + C_\mathrm{HP},
    \end{equation}
    with the length of the heat network \(l_\mathrm{HN}\)
    and the cost of the central heat pump \(C_\mathrm{HP}\)
    calculated according to Eq.~\eqref{eq: HP costs}.
\end{subequations}
As mentioned earlier,
we assumed gas boilers to be existing and thus set \(C_\mathrm{GB} = \SI{0}{\sieuro/a}\).
This is despite the fact that replacements might be needed,
and mainly to have a known Pareto-optimal solution for testing the aggregation methods.

\end{appendices}

\bibliographystyle{unsrt}
\bibliography{references.bib}

\end{document}